\pdfoutput=1
\documentclass[11pt]{amsart}
\usepackage[T1]{fontenc}
\usepackage{lmodern}
\usepackage[a4paper,top=3cm,bottom=3cm,left=2.5cm,right=2.5cm]{geometry}
\usepackage{amssymb,mathtools,bm,microtype}
\usepackage[shortlabels]{enumitem}
\usepackage{longtable}
\usepackage[dvipsnames]{xcolor}
\usepackage[colorlinks=true,citecolor=red,linkcolor=blue,urlcolor=blue]{hyperref}
\hypersetup{
  pdftitle={A power lower bound for the Hilbert numbers of planar polynomial vector fields},
  pdfauthor={Jihao Liu, Chaoyang Qin, Xiaoming Sun},
  pdfsubject={Hilbert's sixteenth problem, limit cycles, planar polynomial vector fields}
}
\ifdefined\pdfinfoomitdate\pdfinfoomitdate=1\fi
\ifdefined\pdftrailerid\pdftrailerid{}\fi
\numberwithin{equation}{section}
\theoremstyle{plain}
\newtheorem{thm}{Theorem}[section]
\newtheorem{lem}[thm]{Lemma}
\newtheorem{prop}[thm]{Proposition}
\newtheorem{cor}[thm]{Corollary}

\theoremstyle{definition}
\newtheorem{defn}[thm]{Definition}

\newtheorem{rem}[thm]{Remark}

\newtheorem{setup}[thm]{Set-up}
\newtheorem{cons}[thm]{Construction}

\newcommand{\id}{\operatorname{id}}
\newcommand{\Rr}{\mathbb{R}}
\newcommand{\ZZ}{\mathbb{Z}}
\providecommand{\doi}[1]{\href{https://doi.org/#1}{doi:#1}}

\begin{document}

\title[Power lower bound for Hilbert numbers]{A power lower bound for the Hilbert numbers of planar polynomial vector fields}

\author{Jihao Liu}
\address[J.~Liu]{School of Mathematical Sciences, Peking University, No.~5 Yiheyuan Road, Haidian District, Beijing 100871, China; and Beijing International Center for Mathematical Research, Peking University, No.~5 Yiheyuan Road, Haidian District, Beijing 100871, China}
\email{liujihao@math.pku.edu.cn}

\author{Chaoyang Qin}
\address[C.~Qin]{State Key Laboratory of Processors, Institute of Computing Technology, Chinese Academy of Sciences, Beijing 100190, China; and School of Computer Science and Technology, University of Chinese Academy of Sciences, Beijing 100049, China}
\email{qinchaoyang25s@ict.ac.cn}

\author{Xiaoming Sun}
\address[X.~Sun]{State Key Laboratory of Processors, Institute of Computing Technology, Chinese Academy of Sciences, Beijing 100190, China; and School of Computer Science and Technology, University of Chinese Academy of Sciences, Beijing 100049, China}
\email{sunxiaoming@ict.ac.cn}

\subjclass[2020]{34C07, 34C05, 37C27}
\keywords{Hilbert's sixteenth problem, Hilbert number, limit cycles, planar polynomial vector fields, formal verification}

\begin{abstract}
For every sufficiently large integer $D$ we construct a real planar polynomial vector field of degree at most $D$ with at least a constant multiple of $D^{2.221004}$ hyperbolic limit cycles. In particular, the Hilbert numbers grow at least like a power of the degree with an exponent larger than $2$.
\end{abstract}

\maketitle

\section{Introduction}\label{sec:intro}

The second part of Hilbert's 16th problem \cite{Hil00,Hil02} asks for the maximal number and the position of the limit cycles of a real planar polynomial vector field of given degree. See \cite{Ily02} for its history. For $D\in\mathbb N$ let $H(D)\in\mathbb N\cup\{\infty\}$ be the supremum of the number of limit cycles (isolated periodic orbits) of real planar polynomial vector fields $X=(X_1,X_2)$ with $\max(\deg X_1,\deg X_2)\le D$. The precise definitions are given in Definition~\ref{defn:cycle}. We prove that $H(D)$ grows at least like a power of $D$ with exponent larger than $2$.

\begin{thm}\label{thm:main}
Let $\lambda^*:=101+\sqrt{10593}$, $q:=\sqrt{\lambda^*}$, and $\beta^*:=2\log 367/\log\lambda^*$. Let $D\ge 343$ be an integer. Then there exist $j\in\mathbb N^+$ with
\[
24\,q^{j}\le D<24\,q^{j+1},
\]
a real planar polynomial vector field $X$ of degree at most $D$, and $367^{j}$ pairwise distinct hyperbolic limit cycles of $X$. Moreover
\[
367^{j}>\Big(\frac{D}{343}\Big)^{\beta^*}\qquad\text{and}\qquad \beta^*>2.221004 .
\]
\end{thm}

\begin{cor}\label{cor:asymptotic}
For every integer $D\ge 343$ there exists a real planar polynomial vector field of degree at most $D$ with at least $343^{-2.221004}D^{2.221004}$ pairwise distinct hyperbolic limit cycles. In particular $H(D)\ge 343^{-2.221004}D^{2.221004}$ for every integer $D\ge343$, and $H(D)=\Omega(D^{2.221004})$.
\end{cor}

Here $\beta^*=2.22100475\ldots$. For $D=343$ the inequality $24q^{j}\le D<24q^{j+1}$ forces $j=1$ (since $24q=342.72\ldots$), so the field of Theorem~\ref{thm:main} has at least $367$ hyperbolic limit cycles.

\emph{Previous lower bounds.} Christopher and Lloyd \cite{CL95} proved that $H(D)$ grows at least like $D^2\log D$. Itenberg and Shustin \cite{IS00} obtained the order $D^2\log D$ for every degree, Li, Chan and Chung \cite{LCC02} improved the explicit bounds, and Han and Li \cite{HL12} reached
\[
\liminf_{D\to\infty}\frac{H(D)}{(D+2)^2\log(D+2)}\ge\frac{1}{2\log2}.
\]
Canard cycles of slow-fast systems reach the leading order $D^2\log D$ by an independent mechanism \cite{ACD+20}. Qin and Sun \cite{QS26} proved $H(D)=\Omega(D^2\log^2D)$, and Lu \cite{Lu26} improved the constant to $H(D)\ge D^2(\log D)^2/(8(\log2)^2)-O(D^2\log D)$. Eshkobilov, Kadyrov and Mamayusupov \cite{EKM26} showed that repeated separable Chebyshev pullbacks, which replicate cycles without creating new ones, yield at most quadratically many cycles in the final degree. We have not found in the literature a lower bound of order $D^{2+\epsilon}$ for any $\epsilon>0$. Corollary~\ref{cor:asymptotic} improves the growth exponent, not the counts in any practical range: the bound $(D/343)^{\beta^*}$ stays below the leading term $D^2(\log D)^2/(8(\log2)^2)$ of \cite{Lu26} for $\log D\le 93.6$.

Individual finiteness was proved by Ilyashenko \cite{Ily91} and \'Ecalle \cite{Eca92} (see also \cite{Liu26}). A proof that $H(D)$ is finite for every $D$ has been announced in \cite{OAI26}. Lower bounds for $H(D)$ do not depend on the finiteness of $H(D)$.

\emph{Idea of the proof.} Up to the last step, the construction of the field in Theorem~\ref{thm:main} is one-dimensional. A segment is a map $x\rightarrow a+cx^{\nu}$ with $c,\nu>0$. A word is a composition of segments. At depth $j$ we obtain a word $P$ with exponent product $1$, length
\[
|P|<5.07\,(\lambda^*)^{j},
\]
and $367^{j}$ fixed points $p$ with $P'(p)\neq1$ (see Theorem~\ref{thm:compiler}, Proposition~\ref{prop:terminal}, Proposition~\ref{prop:iteration}, and Lemma~\ref{lem:perron}).

The iteration acts on analytic one-parameter families of words. One step composes the family with explicit power flows conjugated by members of the family, so that the first variation $V$ of the new family satisfies $V\circ\chi=\rho\,Q$, where $\chi$ is an explicit strictly increasing change of coordinate, $\rho>0$, and $Q$ is $C^1$-close to $\psi\circ(\Lambda V_0/B)$. Here $V_0$ is the first variation of the old family, $B$ the base word of the old family, $\Lambda>0$ a constant, and $\psi$ a small deformation of a translate, multiplied by a positive function, of an explicit $49$-term exponential sum $H_{49}$ with $367$ simple zeros (see Lemma~\ref{lem:phase}, Lemma~\ref{lem:chainvar}, Proposition~\ref{prop:approx}, and Theorem~\ref{thm:cert}). Thus each simple zero of $V_0$ splits into $367$ simple zeros of $V$, while the length budget grows by a linear recursion with Perron root $\lambda^*$. Unlike the replication schemes of \cite{EKM26}, every step creates new simple zeros.

A planar host realizes any word $P$ with exponent product $1$ and $|P|\le4s^2$ as the $C^1$-limit of rescaled first-return maps of polynomial fields of degree at most $20s$ along a polycycle with $4s^2+2$ hyperbolic saddles. One such field turns each fixed point $p$ with $P'(p)\ne1$ into a hyperbolic limit cycle (see Theorem~\ref{thm:host}). With $s\approx\frac12\sqrt{|P|}$ the degree is at most $24q^{j}$ (see Lemma~\ref{lem:arith}) while the number of cycles is $367^{j}=(q^{j})^{\beta^*}$, which yields the exponent $\beta^*$.

\emph{Structure.} Section~\ref{sec:frame} contains the definitions, states the three main ingredients and proves Theorem~\ref{thm:main} and Corollary~\ref{cor:asymptotic}. Section~\ref{sec:cert} proves Theorem~\ref{thm:cert}, Sections~\ref{sec:compilerA} and~\ref{sec:compilerB} prove Theorem~\ref{thm:compiler} and Proposition~\ref{prop:terminal}, and Sections~\ref{sec:hostA} and~\ref{sec:hostB} prove Theorem~\ref{thm:host}. Appendix~\ref{app:data} lists the data of Theorem~\ref{thm:cert}, and Appendix~\ref{app:lean} relates the numbered statements to the formalization.

\emph{Formalization.} A Lean~4 formalization based on Mathlib accompanies the paper and is available at \url{https://github.com/AscendFoam/hilbert16-power-lower-bound-lean}. Its definitions of periodic solution, periodic orbit, isolated orbit and hyperbolic limit cycle, in the file \texttt{Hilbert16/Target.lean}, are those of Definition~\ref{defn:cycle}. The theorem \texttt{hilbert16\_\allowbreak power\_\allowbreak lower\_\allowbreak bound} in \texttt{Hilbert16/Main.lean} states that there exist $c>0$ and $D_0\in\mathbb N$ such that for every integer $D\ge D_0$ there exist a field $X$ of degree at most $D$ and $k\ge cD^{2.221004}$ pairwise distinct hyperbolic limit cycles of $X$. Its proof takes $c=343^{-2.221004}$ and $D_0=343$, as in Corollary~\ref{cor:asymptotic}. The theorem \texttt{finite\_witnesses} in \texttt{Hilbert16/Paper/FiniteWitnesses.lean} is Theorem~\ref{thm:main} without the inequality $\beta^*>2.221004$, which is \texttt{exponent\_lower\_bound} in \texttt{Hilbert16/Paper/Exponent.lean}. The formalization contains the integer computation of Section~\ref{sec:verification}, in exact integer arithmetic, together with a proof that its output encloses the values of $H_{49}$ and $H_{49}'$. Appendix~\ref{app:lean} maps each numbered definition, construction and result to declarations, except Construction~\ref{cons:path}, which has no single counterpart.

\subsection*{Author contributions}
Chaoyang Qin and Xiaoming Sun developed the mathematical construction and proofs. Jihao Liu carried out the Lean~4 formalization and wrote the manuscript. All authors reviewed and approved the final version.

\subsection*{AI-assistance disclosure}
GPT-6.1-sol assisted in exploring and developing the main mathematical results and arguments of this paper, and the SHEAF system assisted in producing the accompanying Lean~4 formalization. The authors take full responsibility for the mathematical correctness and the final text of the manuscript.

\subsection*{Acknowledgements}
The authors thank Ruochuan Liu for his constant support and encouragement. The work was partially supported by the National Key R\&D Program of China \#2024YFA1014400.
\section{The reduction to three theorems}\label{sec:frame}

\subsection{Vector fields and hyperbolic limit cycles}\label{sec:defs}

\begin{defn}\label{defn:cycle}
Let $X=(X_1,X_2)$ be a pair of real polynomials in two variables, viewed as the vector field $z\to X(z)$ on $\Rr^2$, and let $\deg X:=\max(\deg X_1,\deg X_2)$ (total degrees).
\begin{enumerate}
\item A \emph{periodic solution} of $X$ with period $\vartheta>0$ is a differentiable map $\gamma:\Rr\to\Rr^2$ with $\gamma'(t)=X(\gamma(t))$ and $\gamma(t+\vartheta)=\gamma(t)$ for all $t\in\Rr$, which is not constant.
\item A \emph{periodic orbit} of $X$ is a set $\gamma(\Rr)$ with $\gamma$ a periodic solution.
\item A set $O\subset\Rr^2$ is \emph{isolated} if there exists an open set $U\supset O$ such that every periodic orbit $O'$ of $X$ with $O'\subset U$ satisfies $O'=O$.
\item A \emph{hyperbolic limit cycle} of $X$ is an isolated set $O$ for which there exist a periodic solution $\gamma$ with some period $\vartheta$ and $\gamma(\Rr)=O$ such that $\int_0^{\vartheta}\operatorname{div}X(\gamma(t))\,dt\neq0$.
\end{enumerate}
\end{defn}

\begin{rem}\label{rem:cycle}
The period $\vartheta$ in Definition~\ref{defn:cycle}(4) need not be minimal. If $\vartheta$ is $m$ times the minimal period, the integral is $m$ times the integral over the minimal period, so the condition does not depend on $\vartheta$. Since the derivative of the Poincar\'e map of $O$ equals $\exp\int_0^{\vartheta}\operatorname{div}X(\gamma(t))\,dt$ for $\vartheta$ the minimal period \cite[Section~3.4, Theorem~2]{Per01}, the condition says that the derivative of the Poincar\'e map of $O$ is $\neq1$.
\end{rem}

\subsection{Positive-power words}\label{sec:words}

\begin{defn}\label{defn:word}
A \emph{segment} is a map $g(x)=a+cx^{\nu}$ with $a\in\Rr$, $c>0$, $\nu>0$. A \emph{word} is a finite list $W=(g_1,\dots,g_n)$ of segments, of length $|W|:=n\ge0$. It acts as $W:=g_n\circ\dots\circ g_1$, with prefixes $W^{(k)}:=g_k\circ\dots\circ g_1$ ($W^{(0)}=\id$) and exponent product $\Pi(W):=\nu_1\cdots\nu_n$ ($\Pi(\emptyset)=1$). For words $W,Z$, $W\ast Z$ is the list $W$ followed by $Z$. It acts as $Z\circ W$ and $\Pi(W\ast Z)=\Pi(W)\Pi(Z)$. A word $W$ is \emph{legal} on $K\subset\Rr$ if $W^{(k)}(x)>0$ for all $x\in K$ and $0\le k\le n$.
\end{defn}

\begin{lem}\label{lem:legal}
Let $W$ be a word.
\begin{enumerate}
\item The set $\mathcal D_W:=\{x: W\text{ is legal on }\{x\}\}$ is open and contains $(R,\infty)$ for some $R>0$, and $W$ is real analytic with $W>0$ and $W'>0$ on $\mathcal D_W$.
\item Let $U$ be a word, $O\subset(0,\infty)$ a nonempty open interval, and assume that $U$ is legal on $O$, $W$ is legal on $U(O)$ and $W(U(x))=x$ for $x\in O$. Then $\Pi(W)\Pi(U)=1$, and $W\circ U=\id$ on $\mathcal D_{U\ast W}$.
\end{enumerate}
\end{lem}

\begin{proof}
(1) Let $W=(g_1)\ast W_1$. Then $\mathcal D_W=\{x>0: g_1(x)\in\mathcal D_{W_1}\}$, and $g_1$ is real analytic with $g_1'>0$ on $(0,\infty)$ and $g_1(x)\to\infty$ as $x\to\infty$. Induction on $|W|$, starting from $\mathcal D_\emptyset=(0,\infty)$, shows (1) together with: $\mathcal D_W$ is an upper set, hence an open half-line.

(2) The word $U\ast W$ is legal on $O$, so $O\subset\mathcal D_{U\ast W}$. By (1), $W\circ U$ is real analytic on the half-line $\mathcal D_{U\ast W}$ and equals $\id$ on $O$. By the identity principle $W\circ U=\id$ on $\mathcal D_{U\ast W}\supset(R,\infty)$. Induction on the length shows that for every word $Z$ there exists $A_Z>0$ with $Z(x)/x^{\Pi(Z)}\to A_Z$ as $x\to\infty$. For $Z=U\ast W$ the ratio $Z(x)/x^{\Pi(Z)}$ is $x^{1-\Pi(U)\Pi(W)}$ for $x>R$, so $\Pi(U)\Pi(W)=1$.
\end{proof}

\begin{lem}\label{lem:stability}
Let $x_-<x_+$, let $f$ be $C^1$ on an open set containing $[x_-,x_+]$, and let $p_1,\dots,p_K\in(x_-,x_+)$ be pairwise distinct with $f(p_i)=0\neq f'(p_i)$. Then there exist $\epsilon>0$ and pairwise disjoint intervals $[l_i,u_i]\subset(x_-,x_+)$ with $l_i<p_i<u_i$ such that every $g$ that is $C^1$ on an open set containing $[x_-,x_+]$ with $\sup_{[x_-,x_+]}|g-f|<\epsilon$ and $\sup_{[x_-,x_+]}|g'-f'|<\epsilon$ has in each $[l_i,u_i]$ exactly one zero $q_i$. Moreover $q_i\in(l_i,u_i)$ and $g'(q_i)\ne0$.
\end{lem}

\begin{proof}
Let $c_i:=|f'(p_i)|$. By continuity of $f'$ there exist pairwise disjoint $[l_i,u_i]\subset(x_-,x_+)$ with $l_i<p_i<u_i$ and $|f'-f'(p_i)|<c_i/2$ on $[l_i,u_i]$. Then $f(l_i)f(u_i)<0$. Let
\[
\epsilon:=\tfrac12\min_{1\le i\le K}\min\big(|f(l_i)|,|f(u_i)|,c_i/2\big).
\]
For $g$ as in the lemma, $g(l_i)g(u_i)<0$, and $|g'-f'(p_i)|<c_i/2+c_i/4$ on $[l_i,u_i]$, so $g'$ has the strict sign of $f'(p_i)$ there. The intermediate value theorem and the mean value theorem \cite[Theorems~4.23 and 5.10]{Rud76} imply the lemma.
\end{proof}

\begin{defn}\label{defn:family}
Let $n,b,M\in\mathbb N$. An element of $\mathcal F(n,b,M)$ consists of the following data.
\begin{enumerate}
\item Reals $0<x_-\le x_+$, $I:=[x_-,x_+]$, and $\xi_0>0$.
\item Word families $\xi\to W_\xi$ and $\xi\to W^{-}_\xi$ of lengths at most $n$ whose segment coefficients are real analytic functions of $\xi$ on a neighbourhood of $[-\xi_0,\xi_0]$, such that for $|\xi|\le\xi_0$: $W_\xi$ is legal on $I$, $W^{-}_\xi$ is legal on $W_\xi(I)$, $W^{-}_\xi(W_\xi(x))=x$ for $x\in I$, and $\Pi(W_\xi)=1$.
\item Words $B,B^{-}$ of lengths at most $b$ with $B=W_0$ on $I$, $B$ legal on $I$, $B^{-}$ legal on $B(I)$, $B^{-}(B(x))=x$ for $x\in I$, and $\Pi(B)=1$.
\item Points $x_1<\dots<x_M$ in $(x_-,x_+)$ with $V(x_r)=0\ne V'(x_r)$, where $V(x):=\partial_\xi W_\xi(x)|_{\xi=0}$.
\end{enumerate}
\end{defn}

If $M\ge1$, then $x_-<x_1<x_+$. If $n\le n'$ and $b\le b'$, then $\mathcal F(n,b,M)\subset\mathcal F(n',b',M)$.

\subsection{The three main ingredients}\label{sec:modules}

For $N\in\mathbb N^+$, $a,\gamma\in\Rr$ and $\alpha,\beta\in\Rr^N$ let
\begin{equation}\label{eq:template}
H_N(z):=\gamma\big(e^{az}-e^{-az}\big)+\sum_{i=1}^{N}\alpha_i\big(e^{az-\beta_ie^{z}}-e^{-az-\beta_ie^{-z}}\big)\qquad(z\in\Rr).
\end{equation}

\begin{thm}\label{thm:cert}
Let $N=49$ and let $a,\gamma,\alpha_i,\beta_i$ $(1\le i\le 49)$ be the dyadic rationals given in Appendix~\ref{app:data}. Then $a,\gamma,\alpha_i,\beta_i>0$, $\gamma=2^{-360}$, and $H_{49}$ has $367$ pairwise distinct zeros $z_1,\dots,z_{367}\in(-22,22)$ with $H_{49}'(z_k)\ne0$ for $1\le k\le367$.
\end{thm}

\begin{thm}\label{thm:compiler}
Let $n,b,M,R,N\in\mathbb N$ with $N\ge1$, let $W\in\mathcal F(n,b,M)$, let $a>0$, $\gamma\in\Rr$, $\alpha\in\Rr^N$, and $\beta\in(0,\infty)^N$, and let $z_1,\dots,z_R$ be pairwise distinct reals with $H_N(z_k)=0\neq H_N'(z_k)$ for $1\le k\le R$. Then there exist $n',b'\in\mathbb N$ with
\[
n'\le(4N+4)n+(4N+2)b+16N+12,\qquad b'\le 4n+2b+8,
\]
such that $\mathcal F(n',b',RM)$ is nonempty.
\end{thm}

\begin{prop}\label{prop:terminal}
Let $n,b\in\mathbb N$, $M\in\mathbb N^+$, and $W\in\mathcal F(n,b,M)$, with the notation of Definition~\ref{defn:family}. Then there exists $u>0$ such that the word $P:=B^{-}\ast W_u$ (acting as $W_u\circ B^{-}$) satisfies:
\begin{enumerate}
\item $0<B(x_-)<B(x_+)$, $P$ is legal on $[B(x_-),B(x_+)]$, $|P|\le n+b$, and $\Pi(P)=1$.
\item There exist pairwise distinct $p_1,\dots,p_M\in(B(x_-),B(x_+))$ with $P(p_r)=p_r$ and $P'(p_r)\neq1$.
\end{enumerate}
\end{prop}

\begin{thm}\label{thm:host}
Let $0<x_-\le x_+$, let $W$ be a word legal on $[x_-,x_+]$ with $\Pi(W)=1$, and let $s:=\max\big(2,\lceil\sqrt{|W|}/2\rceil\big)$. Let $M\in\mathbb N$ and let $p_1,\dots,p_M\in(x_-,x_+)$ be pairwise distinct with $W(p_r)=p_r$ and $W'(p_r)\neq1$. Then there exist a real planar polynomial vector field $X$ with $\deg X\le 20s$ and $M$ pairwise distinct hyperbolic limit cycles $O_1,\dots,O_M$ of $X$.
\end{thm}

Theorem~\ref{thm:cert} is proved in Section~\ref{sec:cert}, Theorem~\ref{thm:compiler} in Sections~\ref{sec:compilerA}--\ref{sec:compilerB}, Proposition~\ref{prop:terminal} in Section~\ref{sec:terminal}, and Theorem~\ref{thm:host} in Sections~\ref{sec:hostA}--\ref{sec:hostB}.

\subsection{Iteration and arithmetic}\label{sec:arith}

\begin{lem}\label{lem:seed}
Let $I=[1,3]$, $\xi_0=\log(4/3)$, let $W_\xi$ be the one-segment word $x\to(2-2e^{\xi})+e^{\xi}x$, $W^{-}_\xi:=W_{-\xi}$, $B=B^{-}=\emptyset$, and $x_1=2$. The data $(I,\xi_0,W_\xi,W^{-}_\xi,B,B^{-},x_1)$ form an element of $\mathcal F(1,0,1)$.
\end{lem}

\begin{proof}
The coefficients $2-2e^{\xi}$, $e^{\xi}>0$, $\nu=1$ are real analytic on $\Rr$, and $\Pi(W_\xi)=1$. For $|\xi|\le\xi_0$ and $x\in I$,
\[
W_\xi(x)=2+e^{\xi}(x-2)\ge2-\tfrac43>0,\qquad W_{-\xi}(W_\xi(x))=x\ge1 ,
\]
so $W_\xi$ is legal on $I$ and $W^{-}_\xi$ on $W_\xi(I)$. Since $W_0=\id$ and $I\subset(0,\infty)$, the empty words satisfy Definition~\ref{defn:family}(3). Finally $V(x)=x-2$, so $V(2)=0\neq1=V'(2)$.
\end{proof}

\begin{prop}\label{prop:iteration}
Let $(n_r,b_r)_{r\in\mathbb N}$ be defined by $(n_0,b_0)=(1,0)$ and
\[
n_{r+1}=200\,n_r+198\,b_r+796,\qquad b_{r+1}=4\,n_r+2\,b_r+8 .
\]
Then for every $j\in\mathbb N$ there exist reals $0<x_-<x_+$, a word $P$ legal on $[x_-,x_+]$ with $1\le|P|\le n_j+b_j$ and $\Pi(P)=1$, and $367^{j}$ pairwise distinct points $p_r\in(x_-,x_+)$ with $P(p_r)=p_r$ and $P'(p_r)\ne1$.
\end{prop}

\begin{proof}
We show $\mathcal F(n_r,b_r,367^r)\neq\emptyset$ by induction on $r$. The case $r=0$ is Lemma~\ref{lem:seed}. For the step, Theorem~\ref{thm:cert} provides $a>0$, $\gamma$, $\alpha$, $\beta\in(0,\infty)^{49}$, and $367$ simple zeros of $H_{49}$, and Theorem~\ref{thm:compiler} with $N=49$, $R=367$, $M=367^r$ yields $\mathcal F(n',b',367^{r+1})\neq\emptyset$ with
\[
n'\le 200\,n_r+198\,b_r+796=n_{r+1},\qquad b'\le 4n_r+2b_r+8=b_{r+1},
\]
since $(4\cdot49+4,\,4\cdot49+2,\,16\cdot49+12)=(200,198,796)$. Enlarging the length bounds implies $\mathcal F(n_{r+1},b_{r+1},367^{r+1})\ne\emptyset$.

For an element of $\mathcal F(n_j,b_j,367^j)$ with interval $[y_-,y_+]$ and base word $B$, Proposition~\ref{prop:terminal} with $M=367^j\ge1$ yields $P$ and the points $p_r$ on $[x_-,x_+]:=[B(y_-),B(y_+)]$ with $|P|\le n_j+b_j$. If $|P|=0$, then $P=\id$ and $P'(p_1)=1$. Thus $|P|\ge1$.
\end{proof}

\begin{lem}\label{lem:perron}
Let $(n_r)_{r\in\mathbb N},(b_r)_{r\in\mathbb N}$ be nonnegative reals with $n_0=1$, $b_0=0$, $n_{r+1}\le 200n_r+198b_r+796$, and $b_{r+1}\le 4n_r+2b_r+8$ for all $r$. Then $n_j+b_j<\frac{507}{100}(\lambda^*)^{j}$ for every $j\in\mathbb N^+$.
\end{lem}

\begin{proof}
From $102922^2<10593\cdot10^6<102923^2$ we get $102.922<\sqrt{10593}<102.923$. Let $\eta^*:=(\sqrt{10593}-99)/4$, so $49/50<\eta^*<1$, and $C:=796+8\eta^*<804$. A direct calculation shows
\[
200+4\eta^*=\lambda^*,\qquad 198+2\eta^*=\lambda^*\eta^* .
\]
Hence $Z_r:=n_r+\eta^*b_r$ satisfies $Z_0=1$ and $Z_{r+1}\le\lambda^*Z_r+C$, so
\[
Z_j\le(\lambda^*)^j+C\,\frac{(\lambda^*)^j-1}{\lambda^*-1}<\Big(1+\frac{C}{\lambda^*-1}\Big)(\lambda^*)^j .
\]
Since $n_j\ge0$ and $\eta^*<1$, $\eta^*(n_j+b_j)\le Z_j$. With $\lambda^*-1>5073/25$,
\[
\frac{1}{\eta^*}\Big(1+\frac{C}{\lambda^*-1}\Big)<\frac{50}{49}\Big(1+\frac{20100}{5073}\Big)=\frac{1258650}{248577}<\frac{507}{100},
\]
the last inequality being $125865000<126028539$. The lemma follows.
\end{proof}

\begin{lem}\label{lem:arith}
Let $\lambda^*,q,\beta^*$ be as in Theorem~\ref{thm:main}.
\begin{enumerate}
\item $q>14$, $24q<343$, and $(q^{j})^{\beta^*}=367^{j}$ for every $j\in\mathbb N$.
\item For every $j\in\mathbb N^+$ and every real $x$ with $1\le x\le\frac{507}{100}(\lambda^*)^{j}$,
\[
20\max\big(2,\lceil\sqrt{x}/2\rceil\big)<24\,q^{j}.
\]
\item For every integer $D\ge343$ there exists $j\in\mathbb N^+$ with $24q^{j}\le D<24q^{j+1}$.
\item For every real $D>0$ and $j\in\mathbb N$ with $D<24q^{j+1}$, $(D/343)^{\beta^*}<367^{j}$.
\end{enumerate}
\end{lem}

\begin{proof}
(1) From $100<\sqrt{10593}<102.923$: $q^2=\lambda^*>196$ and $(24q)^2=576\lambda^*<117459.648<343^2$. Since $\log q=\frac12\log\lambda^*$, $\beta^*\log q=\log367$, that is $q^{\beta^*}=367$. Hence $(q^j)^{\beta^*}=(q^{\beta^*})^j=367^j$.

(2) Let $Q:=q^j\ge q>14$, so $(\lambda^*)^j=Q^2$. Since $(563/250)^2>507/100$, $\sqrt x<\frac{563}{250}Q$. Since $Q>14$, $20<\frac{10}{7}Q$. Hence
\[
20\lceil\sqrt x/2\rceil<10\sqrt x+20<\tfrac{4191}{175}\,Q<24\,Q,
\]
and $40<24\cdot14<24\,Q$. Thus $20\max(2,\lceil\sqrt x/2\rceil)<24\,Q$.

(3) Since $D/24\ge1$ and $q>1$, there exists $j\in\mathbb N$ with $q^j\le D/24<q^{j+1}$. The case $j=0$ would imply $D<24q<343$.

(4) By (1), $0<D<24q\cdot q^j<343\,q^j$. Since $\beta^*>0$, $t\to t^{\beta^*}$ is strictly increasing on $[0,\infty)$, so $(D/343)^{\beta^*}<(q^j)^{\beta^*}=367^j$.
\end{proof}

\begin{lem}\label{lem:exponent}
$\beta^*>\frac{555251}{250000}=2.221004$.
\end{lem}

\begin{proof}
For $0\le t<1$ let $\Sigma(t):=2\sum_{i=0}^{11}t^{2i+1}/(2i+1)$ and $E(t):=2t^{25}/(1-t^2)$. The series $\log\frac{1+t}{1-t}=2\sum_{i\ge0}t^{2i+1}/(2i+1)$ has nonnegative terms and the tail after $12$ terms is at most $E(t)$, so
\[
\Sigma(t)\le\log\tfrac{1+t}{1-t}\le\Sigma(t)+E(t).
\]
Let $t_1:=1/3$, $t_2:=111/623$, $t_3:=75922301/331922301$, and $U:=203922301/10^6$. Then $\frac{1+t_1}{1-t_1}=2$, $\frac{1+t_2}{1-t_2}=\frac{367}{256}$, and $\frac{1+t_3}{1-t_3}=\frac{U}{128}$. From $102922301^2>10593\cdot10^{12}$ we get $\lambda^*<U$. Thus
\[
\log367\ge8\,\Sigma(t_1)+\Sigma(t_2),\qquad \log\lambda^*<\log U\le7\big(\Sigma+E\big)(t_1)+\big(\Sigma+E\big)(t_3).
\]
Exact rational arithmetic shows $\Sigma(t_1)\ge0.693147$, $\Sigma(t_2)\ge0.360184$, $(\Sigma+E)(t_3)\le0.465710$, and $E(t_1)\le10^{-9}$. Since $8\cdot500000-7\cdot555251=113243$, $7\cdot555251=3886757$ and
\[
113243\cdot693147+500000\cdot360184-555251\cdot465710=102511,
\]
we obtain
\begin{align*}
&500000\log367-555251\log\lambda^*\\
&\qquad>113243\,\Sigma(t_1)+500000\,\Sigma(t_2)-555251\,(\Sigma+E)(t_3)-3886757\,E(t_1)\\
&\qquad\ge0.102511-0.003886757>0.
\end{align*}
Dividing by $250000\log\lambda^*>0$ proves the lemma.
\end{proof}

\begin{proof}[Proof of Theorem~\ref{thm:main}]
Lemma~\ref{lem:arith}(3) yields $j\in\mathbb N^+$ with $24q^j\le D<24q^{j+1}$. Proposition~\ref{prop:iteration} yields $x_-<x_+$, a word $P$ and $367^j$ points $p_r$. The budgets of Proposition~\ref{prop:iteration} satisfy the hypotheses of Lemma~\ref{lem:perron} with equality, so
\[
1\le|P|\le n_j+b_j<\tfrac{507}{100}(\lambda^*)^j .
\]
Theorem~\ref{thm:host} with $W=P$ and $M=367^j$ yields $X$ and $367^j$ pairwise distinct hyperbolic limit cycles with
\[
\deg X\le20\max\big(2,\lceil\sqrt{|P|}/2\rceil\big)<24q^j\le D,
\]
by Lemma~\ref{lem:arith}(2) with $x=|P|$. Lemma~\ref{lem:arith}(4) implies $367^j>(D/343)^{\beta^*}$, and Lemma~\ref{lem:exponent} is $\beta^*>2.221004$. The theorem follows.
\end{proof}

\begin{proof}[Proof of Corollary~\ref{cor:asymptotic}]
Let $D\ge343$ and let $j,X$ be as in Theorem~\ref{thm:main}. Since $D/343\ge1$ and $2.221004<\beta^*$,
\[
343^{-2.221004}D^{2.221004}=\Big(\frac{D}{343}\Big)^{2.221004}\le\Big(\frac{D}{343}\Big)^{\beta^*}<367^j .
\]
A hyperbolic limit cycle is an isolated periodic orbit, so $H(D)\ge367^j$. The corollary follows.
\end{proof}
\section{The 367 simple zeros of the template}\label{sec:cert}

\subsection{Reduction to finitely many inequalities}

\begin{lem}\label{lem:factor}
Let $N\in\mathbb N^+$, $a,\gamma\in\Rr$, $\alpha,\beta\in\Rr^N$, and
\[
\Phi(z):=e^{az}\big(\gamma+\sum_{i=1}^N\alpha_ie^{-\beta_ie^{z}}\big).
\]
Then $H_N$ is odd and real analytic, and for all $z\in\Rr$
\begin{gather*}
H_N(z)=\Phi(z)-\Phi(-z),\qquad H_N'(z)=\Phi'(z)+\Phi'(-z),\\
\Phi'(z)=e^{az}\Big(\gamma a+\sum_{i=1}^N\alpha_ie^{-\beta_ie^{z}}(a-\beta_ie^{z})\Big).
\end{gather*}
\end{lem}

\begin{proof}
Since $e^xe^y=e^{x+y}$, we have
\[
\Phi(w)=\gamma e^{aw}+\sum_{i=1}^N\alpha_ie^{aw-\beta_ie^{w}}.
\]
Comparing with \eqref{eq:template} at $w=z$ and $w=-z$ yields $H_N(z)=\Phi(z)-\Phi(-z)$, hence $H_N(-z)=-H_N(z)$. Every term of \eqref{eq:template} is a constant multiple of $\exp$ composed with a real analytic function, so $H_N$ is real analytic. The product and chain rules yield the formula for $\Phi'$, and $\frac{d}{dz}\Phi(-z)=-\Phi'(-z)$ yields the formula for $H_N'$.
\end{proof}

\begin{lem}\label{lem:oddbracket}
Let $K\in\mathbb N$, $\varrho>0$, and let $f:\Rr\to\Rr$ be differentiable and odd with $f'(0)\ne0$. Let $0<u_1<v_1<u_2<\dots<u_K<v_K<\varrho$ with $f(u_k)f(v_k)<0$ and $f'\neq0$ on $(u_k,v_k)$ for $1\le k\le K$. Then $f$ has $2K+1$ pairwise distinct zeros in $(-\varrho,\varrho)$ at which $f'\ne0$. They are $0$, one point $z_k\in(u_k,v_k)$ for each $k$, and $-z_1,\dots,-z_K$.
\end{lem}

\begin{proof}
Let $1\le k\le K$. Since $f$ is continuous and $f(u_k)f(v_k)<0$, the intermediate value theorem \cite[Theorem~4.23]{Rud76} yields $p_k\in(u_k,v_k)$ with $f(p_k)=0$, and $f'(p_k)\ne0$ by hypothesis. The intervals $(u_k,v_k)$ are pairwise disjoint subsets of $(0,\varrho)$, so $0<p_1<\dots<p_K<\varrho$. Oddness implies $f(0)=0$ and $f(-p_k)=0$, and differentiating $f(-x)=-f(x)$ yields $f'(-x)=f'(x)$ for $x\in\Rr$, so $f'(-p_k)\ne0$. The points $0$, $p_k$, $-p_k$ $(1\le k\le K)$ are pairwise distinct by sign and by $p_1<\dots<p_K$. The lemma follows.
\end{proof}

\begin{prop}\label{prop:checks}
Let $a,\gamma,\alpha_i,\beta_i$ $(1\le i\le49)$ and $u_k,v_k:=u_k+2^{-82}$ $(1\le k\le183)$ be the numbers of Appendix~\ref{app:data}. Then:
\begin{enumerate}
\item $a,\gamma,\alpha_i,\beta_i>0$ and $\gamma=2^{-360}$.
\item $0<u_1<v_1<u_2<\dots<u_{183}<v_{183}<22$.
\item $H_{49}(u_k)\,H_{49}(v_k)<0$ for $1\le k\le183$.
\item For each $k$, either $H_{49}'>0$ on $[u_k,v_k]$ or $H_{49}'<0$ on $[u_k,v_k]$.
\item $H_{49}'(0)<0$.
\end{enumerate}
\end{prop}

\begin{proof}[Proof of Theorem~\ref{thm:cert}]
Positivity and $\gamma=2^{-360}$ are Proposition~\ref{prop:checks}(1). By Lemma~\ref{lem:factor}, $f:=H_{49}$ is differentiable and odd, and $f'(0)\ne0$ by Proposition~\ref{prop:checks}(5). The hypotheses of Lemma~\ref{lem:oddbracket} with $K=183$ and $\varrho=22$ hold: the order of the $u_k,v_k$ by Proposition~\ref{prop:checks}(2), the signs by (3), and $f'\ne0$ on $(u_k,v_k)$ by (4). Thus $H_{49}$ has $2\cdot183+1=367$ pairwise distinct zeros $z_1,\dots,z_{367}\in(-22,22)$ with $H_{49}'(z_k)\ne0$. The theorem follows.
\end{proof}

Proposition~\ref{prop:checks} is proved by a computer-assisted verification in exact integer arithmetic (Section~\ref{sec:verification}).

\subsection{Exact interval arithmetic}\label{sec:interval}

All numbers of Appendix~\ref{app:data} are of the form $m2^{e}$ with $m,e\in\ZZ$ and $e\ge-380$. In particular they are elements of $2^{-960}\ZZ$.

\begin{cons}\label{cons:interval}
A pair
\[
\mathbf X=(X^{\operatorname{lo}},X^{\operatorname{hi}})\in\ZZ^2
\]
denotes
\[
[X^{\operatorname{lo}}2^{-960},X^{\operatorname{hi}}2^{-960}].
\]
It is \emph{ordered} if $X^{\operatorname{lo}}\le X^{\operatorname{hi}}$, and it \emph{encloses} $y\in\Rr$ if $X^{\operatorname{lo}}\le2^{960}y\le X^{\operatorname{hi}}$. For $y\in2^{-960}\ZZ$ let $\langle y\rangle:=(2^{960}y,2^{960}y)$. The operations are as follows.
\begin{enumerate}
\item
\[
\mathbf X+\mathbf Y:=(X^{\operatorname{lo}}+Y^{\operatorname{lo}},X^{\operatorname{hi}}+Y^{\operatorname{hi}}),
\]
$-\mathbf X:=(-X^{\operatorname{hi}},-X^{\operatorname{lo}})$, $\mathbf X-\mathbf Y:=\mathbf X+(-\mathbf Y)$.
\item
\[
\mathbf X\cdot\mathbf Y:=\big(\lfloor2^{-960}\min\mathcal P\rfloor,\lceil2^{-960}\max\mathcal P\rceil\big),
\]
where
\[
\mathcal P:=\{X^{\operatorname{lo}}Y^{\operatorname{lo}},X^{\operatorname{lo}}Y^{\operatorname{hi}},X^{\operatorname{hi}}Y^{\operatorname{lo}},X^{\operatorname{hi}}Y^{\operatorname{hi}}\}.
\]

\item
\[
\operatorname{rec}\mathbf X:=\big(\lfloor2^{1920}/X^{\operatorname{hi}}\rfloor,\lceil2^{1920}/X^{\operatorname{lo}}\rceil\big)
\]
and
\[
\operatorname{sq}\mathbf X:=\big(\lfloor(X^{\operatorname{lo}})^2/2^{960}\rfloor,\lceil(X^{\operatorname{hi}})^2/2^{960}\rceil\big).
\]

\item $\operatorname{sgn}\mathbf X:=1$ if $X^{\operatorname{lo}}>0$, $\operatorname{sgn}\mathbf X:=-1$ if $X^{\operatorname{hi}}<0$. Otherwise $\operatorname{sgn}\mathbf X$ is undefined.
\item For $x\in\ZZ$ let $\mathcal E(0):=\langle1\rangle$, and let $\mathcal E(x):=(0,1)$ if $x\le-960\cdot2^{960}$. Otherwise let $A:=|x|$, let $\ell$ be the least integer with $A<2^{\ell}$, $r:=\max(0,\ell-940)$, $d:=2^{960+r}$,
\[
\tau^{\operatorname{lo}}_0=\tau^{\operatorname{hi}}_0:=2^{960},\qquad \tau^{\operatorname{lo}}_{j+1}:=\Big\lfloor\frac{\tau^{\operatorname{lo}}_jA}{d(j+1)}\Big\rfloor,\qquad \tau^{\operatorname{hi}}_{j+1}:=\Big\lceil\frac{\tau^{\operatorname{hi}}_jA}{d(j+1)}\Big\rceil\qquad(j\in\mathbb N),
\]
let $n$ be the least index with $\tau^{\operatorname{hi}}_n\le1$,
\[
\mathbf T:=\big(\sum_{j=0}^n\tau^{\operatorname{lo}}_j,\,1+\sum_{j=0}^n\tau^{\operatorname{hi}}_j\big),
\]
and
\[
\mathcal E(x):=\operatorname{sq}^r\mathbf T\ \text{ if } x>0,\qquad \mathcal E(x):=\operatorname{sq}^r\operatorname{rec}\mathbf T\ \text{ if } x<0.
\]
\item $\mathcal E(\mathbf X):=\mathcal E(X^{\operatorname{lo}})$ if $X^{\operatorname{lo}}=X^{\operatorname{hi}}$, and
\[
\mathcal E(\mathbf X):=\big(\mathcal E(X^{\operatorname{lo}})^{\operatorname{lo}},\mathcal E(X^{\operatorname{hi}})^{\operatorname{hi}}\big)
\]
otherwise.
\end{enumerate}
\end{cons}

\begin{lem}\label{lem:sound}
Let $\mathbf X,\mathbf Y$ be ordered pairs enclosing $x,y\in\Rr$, and let $m\in\ZZ$. Then:
\begin{enumerate}
\item $\mathbf X+\mathbf Y$, $-\mathbf X$, $\mathbf X-\mathbf Y$, $\mathbf X\cdot\mathbf Y$ are ordered and enclose $x+y$, $-x$, $x-y$, $xy$.
\item If $X^{\operatorname{lo}}>0$, then $\operatorname{rec}\mathbf X$ is ordered, has lower endpoint $\ge0$ and encloses $1/x$. If $X^{\operatorname{lo}}\ge0$, then $\operatorname{sq}\mathbf X$ is ordered, has lower endpoint $\ge0$ and encloses $x^2$.
\item $\mathcal E(m)$ and $\mathcal E(\mathbf X)$ are ordered, have lower endpoints $\ge0$, and enclose $e^{m2^{-960}}$ and $e^{x}$, respectively. If $0\ne m>-960\cdot2^{960}$, then the index $n$ of Construction~\ref{cons:interval}(5) with $x=m$ exists and $1\le n\le48$.
\item If $\operatorname{sgn}\mathbf X$ is defined, then $x\operatorname{sgn}\mathbf X>0$.
\end{enumerate}
\end{lem}

\begin{proof}
(1), (2), (4) follow from $\min\mathcal P\le2^{1920}xy\le\max\mathcal P$, $1/x$ is decreasing and $x^2$ increasing on $(0,\infty)$, and $\lfloor c\rfloor\le c\le\lceil c\rceil$.

(3) Let $m\ne0$ and $m>-960\cdot2^{960}$, so $A=|m|>0$. Then $t:=A/d\in(0,2^{-20})$, since $A<2^{940}=2^{-20}d$ if $r=0$ and $A<2^{\ell}=2^{-20}d$ if $r>0$. By induction on $j$,
\[
0\le\tau^{\operatorname{lo}}_j\le2^{960}t^j/j!\le\tau^{\operatorname{hi}}_j\qquad\text{and}\qquad 1\le\tau^{\operatorname{hi}}_j.
\]
Moreover $\tau^{\operatorname{hi}}_j\le2^{960-20j}$ for $j\le48$, since
\[
\tau^{\operatorname{hi}}_jA/(d(j+1))<\tau^{\operatorname{hi}}_j2^{-20}
\]
and $2^{960-20(j+1)}\in\ZZ$ for $j+1\le48$. Hence $1\le n\le48$ and $\tau^{\operatorname{hi}}_n=1$. By Taylor's theorem \cite[Theorem~5.15]{Rud76} there exists $c\in(0,t)$ with $e^t-\sum_{j=0}^nt^j/j!=e^ct^{n+1}/(n+1)!$, and $e^c<2$, so
\[
0<2^{960}\Big(e^t-\sum_{j=0}^n\frac{t^j}{j!}\Big)<\frac{2t}{n+1}\cdot\frac{2^{960}t^n}{n!}\le2t\,\tau^{\operatorname{hi}}_n<1 .
\]
Thus $\mathbf T$ is ordered, $T^{\operatorname{lo}}\ge2^{960}$, and $\mathbf T$ encloses $e^t$. By (2), $\operatorname{rec}\mathbf T$ encloses $e^{-t}$, and $r$ squarings enclose $e^{\pm2^rt}=e^{m2^{-960}}$. If $m\le-960\cdot2^{960}$, then $0<e^{m2^{-960}}\le e^{-960}<2^{-960}$ because $e>2$. The case $m=0$ is trivial, and the claim on $\mathcal E(\mathbf X)$ follows since $\exp$ is increasing.
\end{proof}

\begin{lem}\label{lem:enclosure}
Let $a,\gamma,\alpha_i,\beta_i$ $(1\le i\le49)$ be the numbers of Appendix~\ref{app:data}, let $\mathbf Z$ be an ordered pair and $\delta\in\{0,1\}$. Let $\mathbf E:=\mathcal E(\mathbf Z)$, $\mathbf t_i:=\langle\beta_i\rangle\cdot\mathbf E$,
\[
\mathbf Q_0:=\langle\gamma\rangle\cdot\langle a\rangle^{\delta}
\]
(the factor $\langle a\rangle^0$ omitted), and for $1\le i\le49$
\[
\mathbf Q_i:=\mathbf Q_{i-1}+\big(\langle\alpha_i\rangle\cdot\mathcal E(-\mathbf t_i)\big)\cdot(\langle a\rangle-\mathbf t_i)^{\delta},\qquad \mathbf F_\delta(\mathbf Z):=\mathbf Q_{49}\cdot\mathcal E(\langle a\rangle\cdot\mathbf Z),
\]
and
\[
\mathbf H_0(\mathbf Z):=\mathbf F_0(\mathbf Z)-\mathbf F_0(-\mathbf Z),\qquad \mathbf H_1(\mathbf Z):=\mathbf F_1(\mathbf Z)+\mathbf F_1(-\mathbf Z).
\]
Then $\mathbf H_0(\mathbf Z)$ and $\mathbf H_1(\mathbf Z)$ are ordered and enclose $H_{49}(z)$ and $H_{49}'(z)$, respectively, for every $z$ enclosed by $\mathbf Z$.
\end{lem}

\begin{proof}
Let $z$ be enclosed by $\mathbf Z$. By Lemma~\ref{lem:sound}, used at each operation, $\mathbf E$, $\mathbf t_i$, $\mathcal E(-\mathbf t_i)$ enclose $e^z$, $\beta_ie^z$, $e^{-\beta_ie^z}$, and by induction on $i$ the ordered pair $\mathbf Q_{49}$ encloses $\gamma+\sum_i\alpha_ie^{-\beta_ie^z}$ if $\delta=0$ and $\gamma a+\sum_i\alpha_ie^{-\beta_ie^z}(a-\beta_ie^z)$ if $\delta=1$. Thus $\mathbf F_0(\mathbf Z)$ and $\mathbf F_1(\mathbf Z)$ enclose $\Phi(z)$ and $\Phi'(z)$, with $\Phi$ as in Lemma~\ref{lem:factor}. Since $-\mathbf Z$ encloses $-z$, Lemma~\ref{lem:factor} implies the claim. Each step uses only that the operands of that step enclose the corresponding real numbers, so the conclusion holds for all $z$ enclosed by $\mathbf Z$ simultaneously.
\end{proof}

\subsection{The verification}\label{sec:verification}

\begin{proof}[Proof of Proposition~\ref{prop:checks}]
Items (1) and (2) are comparisons of the integers $2^{960}a$, $2^{960}\gamma$, $2^{960}\alpha_i$, $2^{960}\beta_i$, $2^{960}u_k$, $2^{960}v_k$ with $0$, with each other and with $22\cdot2^{960}$. The equality $\gamma=2^{-360}$ is part of the data. By Lemma~\ref{lem:enclosure} and Lemma~\ref{lem:sound}(4), items (3)--(5) follow from the following strict signs, for $1\le k\le183$:
\begin{enumerate}
\item[(3$'$)] $\operatorname{sgn}\mathbf H_0(\langle u_k\rangle)$ and $\operatorname{sgn}\mathbf H_0(\langle v_k\rangle)$ are defined and opposite.
\item[(4$'$)] The sign
\[
\operatorname{sgn}\mathbf H_1\big((2^{960}u_k,2^{960}v_k)\big)
\]
is defined. The pair encloses every $z\in[u_k,v_k]$.
\item[(5$'$)] $\operatorname{sgn}\mathbf H_1(\langle0\rangle)=-1$.
\end{enumerate}
Items (3$'$)--(5$'$) consist of $550$ evaluations of Lemma~\ref{lem:enclosure} followed by sign tests, a finite computation with integers. It was carried out and returns (3$'$)--(5$'$), with the sign in (4$'$) equal to $\operatorname{sgn}\mathbf H_0(\langle v_k\rangle)$. It is replayed in the formalization (Appendix~\ref{app:lean}). The proposition follows.
\end{proof}

The computation is reproduced by the script \texttt{verify\_template.py} (using exact Python integer arithmetic), which reads the data file \texttt{hilbert16\_367\_roots.json} (see Appendix~\ref{app:data}), performs Construction~\ref{cons:interval}, checks the positivity in (1) and the tests (2), (3$'$)--(5$'$), and prints the acceptance. The equality $\gamma=2^{-360}$ in (1) is the entry $(1,-360)$ of the data file. For orientation, $a\approx11.1538$, $u_1\approx0.0904$, $u_{181}\approx16.312$, $u_{182}\approx18.767$, $u_{183}\approx20.183$, and $0.075<u_{k+1}-u_k<0.092$ for $1\le k\le180$.
\section{The compiler I: phase, transport and the three-block path}\label{sec:compilerA}

In Sections~\ref{sec:compilerA} and~\ref{sec:compilerB}, $N,R,M,n,b$, $W\in\mathcal F(n,b,M)$ with the notation of Definition~\ref{defn:family}, $a>0$, $\gamma\in\Rr$, $\alpha\in\Rr^N$, $\beta\in(0,\infty)^N$, and $z_1,\dots,z_R$ are as in Theorem~\ref{thm:compiler}, and $M\ge1$ (the case $RM=0$ is settled in Section~\ref{sec:compilerproof}). Thus $x_-<x_+$. Since segments are strictly increasing, the sets $\mathcal D_B$ and $\mathcal D_{W_0}$ of Lemma~\ref{lem:legal}(1) are open half-lines containing $I$. As $B=W_0$ on $I$, the identity principle for real analytic functions shows that $B=W_0$ on the open interval $U_0:=\mathcal D_B\cap\mathcal D_{W_0}$, and $B>0$, $B'>0$ on $U_0$.

\subsection{Phase replacement and root multiplication}\label{sec:phase}

\begin{lem}\label{lem:phase}
Let $\ell\in\mathbb N^+$ with $|z_k|<\ell$ for $1\le k\le R$. Define $\Phi_+(u):=\gamma+\sum_{i=1}^N\alpha_ie^{\beta_i^{+}e^{u}}$ and $\Phi_-(u):=\gamma+\sum_{i=1}^N\alpha_ie^{\beta_i^{-}e^{-u}}$ with $\beta_i^{+}:=-\beta_ie^{-\ell}$, $\beta_i^{-}:=-\beta_ie^{\ell}$, and for $\theta>0$
\[
\psi_\theta(u):=\Phi_+(u)-\exp\Big(2a\ell-2a\,\frac{e^{\theta u}-1}{\theta}\Big)\,\Phi_-(u).
\]
Then there exists $\theta_1>0$ such that for every $\theta\in(0,\theta_1)$ the function $\psi_\theta$ has $R$ pairwise distinct zeros $y_1,\dots,y_R\in(0,2\ell)$ with $\psi_\theta'(y_k)\neq0$.
\end{lem}

\begin{proof}
Let $\widetilde H(u):=e^{-a(u-\ell)}H_N(u-\ell)$. Expanding \eqref{eq:template} with $e^{u-\ell}=e^{-\ell}e^{u}$ and $e^{-(u-\ell)}=e^{\ell}e^{-u}$ shows
\[
\widetilde H(u)=\Phi_+(u)-e^{2a\ell-2au}\,\Phi_-(u)\qquad(u\in\Rr).
\]
The function $\widetilde H$ is real analytic, $z_k+\ell\in(0,2\ell)$ are zeros of $\widetilde H$, and $\widetilde H'(z_k+\ell)=e^{-az_k}H_N'(z_k)\neq0$. Let $E_\theta(u):=\exp(-2a(e^{\theta u}-1)/\theta)$ and $E_0(u):=e^{-2au}$. Then $E_\theta,E_0\in(0,1]$ for $u\ge0$ and $\psi_\theta-\widetilde H=-e^{2a\ell}(E_\theta-E_0)\Phi_-$. For $0\le\theta u\le1$ Taylor's theorem \cite[Theorem~5.15]{Rud76} implies $|e^{\theta u}-1-\theta u|\le(\theta u)^2$ and $|e^{\theta u}-1|\le2\theta u$. With $E_\theta'=-2ae^{\theta u}E_\theta$ and $|e^s-1|\le2|s|$ for $|s|\le1$ one obtains, for $0<\theta\le\theta_2:=\min(1/(2\ell),1/(8a\ell^2))$ and $u\in[0,2\ell]$,
\[
|E_\theta-E_0|\le c_0\theta,\qquad|E_\theta'-E_0'|\le c_1\theta,\qquad c_0:=16a\ell^2,\ c_1:=8a\ell+32a^2\ell^2 .
\]
Hence $|\psi_\theta-\widetilde H|+|\psi_\theta'-\widetilde H'|\le C\theta$ on $[0,2\ell]$, where $C:=e^{2a\ell}(c_0+c_1)\sup_{[0,2\ell]}(|\Phi_-|+|\Phi_-'|)$. Let $\epsilon$ and the intervals $J_k\subset(0,2\ell)$ be given by Lemma~\ref{lem:stability} for $\widetilde H$ on $[0,2\ell]$ at the points $z_k+\ell$, and let $\theta_1:=\min(\theta_2,\epsilon/(C+1))$. For $\theta\in(0,\theta_1)$, Lemma~\ref{lem:stability} with $g=\psi_\theta$ yields one zero $y_k$ in the interior of $J_k$ with $\psi_\theta'(y_k)\neq0$ for each $k$. The points $y_k$ are pairwise distinct since the intervals are disjoint. The lemma follows.
\end{proof}

Fix from now on $\ell$ as in Lemma~\ref{lem:phase} and $\theta\in(0,\min(1,\theta_1))$, and let $\delta_*:=-2a/\theta$ and $A_*:=e^{2a\ell-\delta_*}$, so that
\begin{equation}\label{eq:Astar}
A_*\,e^{\delta_*e^{\theta u}}=\exp\Big(2a\ell-2a\,\frac{e^{\theta u}-1}{\theta}\Big)\qquad(u\in\Rr).
\end{equation}

\begin{lem}\label{lem:rootmult}
Let $x_-<x_+$, let $B,V$ be $C^1$ on an open set containing $[x_-,x_+]$ with $B>0$ on $[x_-,x_+]$, let $x_1,\dots,x_M\in(x_-,x_+)$ be pairwise distinct with $V(x_r)=0\neq V'(x_r)$, and let $y_1,\dots,y_R>0$ be pairwise distinct with $\psi(y_k)=0\neq\psi'(y_k)$ for a function $\psi$ differentiable at the $y_k$. Then there exists $\Lambda_0>0$ such that for every $\Lambda\ge\Lambda_0$ the function $\psi\circ(\Lambda V/B)$ has $RM$ pairwise distinct zeros $x_{r,k}\in(x_-,x_+)$ with $\Lambda V(x_{r,k})/B(x_{r,k})=y_k$ and $(\psi\circ(\Lambda V/B))'(x_{r,k})\ne0$.
\end{lem}

\begin{proof}
On an open set containing $[x_-,x_+]$ the function $g:=V/B$ is $C^1$, and $g(x_r)=0$, $g'(x_r)=V'(x_r)/B(x_r)\neq0$. Let $\epsilon$ and the disjoint intervals $J_r\ni x_r$ be given by Lemma~\ref{lem:stability} for $g$, and let $\Lambda_0:=1+\max_ky_k/\epsilon$. For $\Lambda\ge\Lambda_0$, the function $g-y_k/\Lambda$ differs from $g$ by a constant in $(0,\epsilon)$, so it has in $J_r$ exactly one zero $x_{r,k}$, and $g'(x_{r,k})\neq0$. Then $\Lambda g(x_{r,k})=y_k$, and the chain rule shows
\[
(\psi\circ(\Lambda g))'(x_{r,k})=\Lambda\psi'(y_k)g'(x_{r,k})\neq0.
\]
Points with different $r$ are in disjoint intervals, and $x_{r,k}\neq x_{r,k'}$ for $k\neq k'$ since $y_k\neq y_{k'}$. The lemma follows.
\end{proof}

The first variation $V$ is real analytic near $I$, since $(\xi,x)\to W_\xi(x)$ is jointly real analytic near $\{0\}\times I$. Lemma~\ref{lem:rootmult} with the base word $B$ (which is $C^1$ and positive on $U_0$), this $V$, the zeros $x_1,\dots,x_M$ of $V$ (Definition~\ref{defn:family}(4)), $\psi=\psi_\theta$ and the zeros $y_k\in(0,2\ell)$ of Lemma~\ref{lem:phase} implies the existence of $\Lambda_0$. Fix $\Lambda:=\Lambda_0$ and let $\widehat\psi:=\psi_\theta\circ(\Lambda V/B)$, which has $RM$ pairwise distinct zeros in $(x_-,x_+)$ with nonzero derivative.

\subsection{Coordinates and flows}\label{sec:flows}

\begin{lem}\label{lem:weight}
Let $m\in\mathbb N^+$ and $L>0$, and define $\chi_L(u):=\big(1+\frac{u^{m+1}}{(m+1)L}\big)^{L}$ for $u>0$.
\begin{enumerate}
\item $\chi_L$ is a word of two segments with $\Pi(\chi_L)=(m+1)L$. It maps $(0,\infty)$ bijectively onto $(1,\infty)$, and $\chi_L^{-1}(z)=\big((m+1)L(z^{1/L}-1)\big)^{1/(m+1)}$ is a word of two segments with exponent product $((m+1)L)^{-1}$.
\item Let $B$ be $C^1$ with $B>0$, $B'>0$ on an open set containing $[x_-,x_+]$, and define $\chi:=\chi_L\circ B$, $w:=B^{-m}/B'$, and $\varpi_L:=\big(1+\frac{B^{m+1}}{(m+1)L}\big)^{-1}$. Then $\chi'/\chi=\varpi_L/w$, and there exists $C>0$, independent of $L$, with $|\varpi_L-1|+|\varpi_L'|\le C/L$ on $[x_-,x_+]$.
\end{enumerate}
\end{lem}

\begin{proof}
(1) The segments of $\chi_L$ are $u\to1+u^{m+1}/((m+1)L)$ and $v\to v^{L}$. Those of $\chi_L^{-1}$ are $z\to-(m+1)L+(m+1)Lz^{1/L}$ and $v\to v^{1/(m+1)}$. A direct calculation shows the remaining claims. (2) Since $\log\chi=L\log\big(1+\frac{B^{m+1}}{(m+1)L}\big)$, $\chi'/\chi=\varpi_LB^mB'=\varpi_L/w$ and
\[
\varpi_L-1=-\frac{B^{m+1}}{(m+1)L+B^{m+1}},\qquad\varpi_L'=-\frac{\varpi_L^2B^mB'}{L};
\]
as $0<\varpi_L\le1$, $C:=1+\sup_{[x_-,x_+]}\big(B^{m+1}/(m+1)+B^mB'\big)$ works.
\end{proof}

\begin{lem}\label{lem:powerflow}
Let $m\in\mathbb N^+$ and $\eta,\delta\in\Rr$, and, where the radicand is positive, define
\[
\Psi^{(m)}_{\eta,\delta}(z):=\big(z^{m+1}+(m+1)\delta\eta\big)^{1/(m+1)}.
\]
Then $\Psi^{(m)}_{\eta,\delta}$ is a word of two segments with exponent product $1$, $\Psi^{(m)}_{0,\delta}=\id$,
\[
\Psi^{(m)}_{\eta,\delta'}\circ\Psi^{(m)}_{\eta,\delta}=\Psi^{(m)}_{\eta,\delta+\delta'},
\]
and, for every differentiable $z(\cdot)$ with $z(0)>0$,
\[
\partial_\eta\Psi^{(m)}_{\eta,\delta}(z(\eta))|_{\eta=0}=z'(0)+\delta z(0)^{-m}.
\]
\end{lem}

\begin{proof}
The segments are $z\to(m+1)\delta\eta+z^{m+1}$ and $v\to v^{1/(m+1)}$. The radicand of $\Psi^{(m)}_{\eta,\delta'}$ at $\Psi^{(m)}_{\eta,\delta}(z)$ is $z^{m+1}+(m+1)(\delta+\delta')\eta$, which implies the group law and $\Psi^{(m)}_{0,\delta}=\id$. For $\varrho(\eta):=z(\eta)^{m+1}+(m+1)\delta\eta$ we have $\varrho(0)=z(0)^{m+1}$, and
\[
\partial_\eta\varrho^{1/(m+1)}|_{\eta=0}=\varrho(0)^{-m/(m+1)}\varrho'(0)/(m+1)=z(0)^{-m}(z(0)^mz'(0)+\delta).
\]

\end{proof}

For $m\in\mathbb N^+$, $\omega\in\{1,\theta,-1\}$ and $\delta,\eta\in\Rr$ let $\xi_\omega:=-\omega\Lambda/m$,
\[
F_{\omega,\eta}(\delta):=W_{\xi_\omega}\ast\Psi^{(m)}_{\eta,-\delta}\ast W^{-}_{\xi_\omega},\qquad w_\omega:=W_{\xi_\omega}^{-m}/\partial_xW_{\xi_\omega};
\]
the word $F_{\omega,\eta}(\delta)$ acts as
\[
W^{-}_{\xi_\omega}\circ\Psi^{(m)}_{\eta,-\delta}\circ W_{\xi_\omega}.
\]
Assume that $|\xi_\omega|\le\xi_0$. Then $|F_{\omega,\eta}(\delta)|\le2n+2$. By Lemma~\ref{lem:legal}(2), $\Pi(W^{-}_{\xi_\omega})=1$, so $\Pi(F_{\omega,\eta}(\delta))=1$, and $W^{-}_{\xi_\omega}\circ W_{\xi_\omega}=\id$ on the open set
\[
O_\omega:=\mathcal D_{W_{\xi_\omega}\ast W^{-}_{\xi_\omega}}\supset I,
\]
so, for $x\in O_\omega$,
\[
(W^{-}_{\xi_\omega})'(W_{\xi_\omega}(x))=1/\partial_xW_{\xi_\omega}(x).
\]
Lemma~\ref{lem:legal}(2) with $W_{\xi_\omega}$ and $W^{-}_{\xi_\omega}$ exchanged, on the interval $W_{\xi_\omega}((x_-,x_+))$, implies $W_{\xi_\omega}\circ W^{-}_{\xi_\omega}=\id$ on $\mathcal D_{W^{-}_{\xi_\omega}\ast W_{\xi_\omega}}$. Lemma~\ref{lem:powerflow} therefore implies $F_{\omega,0}(\delta)=\id$ on $O_\omega$. At every $x$ at which the word $F_{\omega,\eta}(\delta)\ast F_{\omega,\eta}(\delta')$ is legal, the word $F_{\omega,\eta}(\delta+\delta')$ is legal and
\[
F_{\omega,\eta}(\delta')(F_{\omega,\eta}(\delta)(x))=F_{\omega,\eta}(\delta+\delta')(x).
\]
Moreover,
\begin{equation}\label{eq:movinginput}
\partial_\eta\big(F_{\omega,\eta}(\delta)(x(\eta))\big)\big|_{\eta=0}=x'(0)-\delta\,w_\omega(x(0))
\end{equation}
for every differentiable $x(\cdot)$ with $x(0)\in O_\omega$.

\begin{lem}\label{lem:transport}
There exist an open set $U\supset I$, $m_0\in\mathbb N^+$, and $C>0$ such that for every integer $m\ge m_0$ and $\omega\in\{1,\theta,-1\}$: $|\xi_\omega|<\xi_0$, the functions $W_{\xi_\omega}$, $B$ are analytic with positive derivative on $U$, and for $x\in I$
\[
\Big|\log\frac{w_\omega}{w}-\frac{\omega\Lambda V}{B}\Big|+\Big|\Big(\log\frac{w_\omega}{w}\Big)'-\Big(\frac{\omega\Lambda V}{B}\Big)'\Big|\le\frac{C}{m}.
\]
\end{lem}

\begin{proof}
The set of $(\xi,x)$ at which all prefixes of $W_\xi$ are positive is open (the segment coefficients are analytic in $\xi$) and contains $\{0\}\times I$. By the tube lemma there exist $r\in(0,\xi_0)$ and an open $U$ with $I\subset U\subset U_0$ such that $(\xi,x)\to W_\xi(x)$ is real analytic on $(-r,r)\times U$ with $W_\xi>0$ and $\partial_xW_\xi>0$ there. Let $\mathcal L_0(\xi,x):=\log W_\xi(x)$ and $\mathcal L_1(\xi,x):=\log\partial_xW_\xi(x)$. On $U$, $B=W_0$, hence $w=W_0^{-m}/W_0'$, $\partial_\xi\mathcal L_0(0,\cdot)=V/B$, and
\[
\log\frac{w_\omega}{w}=-m\big(\mathcal L_0(\xi_\omega,\cdot)-\mathcal L_0(0,\cdot)\big)-\big(\mathcal L_1(\xi_\omega,\cdot)-\mathcal L_1(0,\cdot)\big).
\]
Taylor's theorem \cite[Theorem~5.15]{Rud76} in $\xi$ for $\mathcal L_0$, $\partial_x\mathcal L_0$, $\mathcal L_1$, and $\partial_x\mathcal L_1$ on the compact set $[-r/2,r/2]\times I$ yields $A\ge0$ with
\[
\big|\partial_x^{\iota}\big(\mathcal L_0(\xi,x)-\mathcal L_0(0,x)-\xi V(x)/B(x)\big)\big|\le A\xi^2,\qquad\big|\partial_x^{\iota}\big(\mathcal L_1(\xi,x)-\mathcal L_1(0,x)\big)\big|\le A|\xi|
\]
for $\iota\in\{0,1\}$, $|\xi|\le r/2$ and $x\in I$. Let $m_0\in\mathbb N^+$ with $m_0\ge2\Lambda/r$. As $0<\theta<1$, $|\xi_\omega|\le\Lambda/m\le r/2$ for $m\ge m_0$, and $-m\xi_\omega V/B=\omega\Lambda V/B$. Substituting $\xi=\xi_\omega$ shows the bound with $C:=2A(\Lambda^2+\Lambda)+1$.
\end{proof}

\subsection{The three-block path}\label{sec:path}

\begin{cons}\label{cons:path}
Let $m\ge m_0$ be an integer (Lemma~\ref{lem:transport}), $L>0$, $\chi:=\chi_L\circ B$, and let $\chi^{-1}$ denote the word $\chi_L^{-1}\ast B^{-}$, which acts as $B^{-}\circ\chi_L^{-1}$ and satisfies $\chi^{-1}\circ\chi=\id$ on $\mathcal D_{B\ast B^-}$ by Lemma~\ref{lem:legal}(2). Let the path data be as follows, with $\beta^{\pm}_i$ as in Lemma~\ref{lem:phase}.
\begin{enumerate}
\item Steps $(\omega_j,\delta_j)$, $1\le j\le2N+2$: $\omega_j=1$ for $j\le N+1$, with $\delta_1=\beta^+_1$, $\delta_j=\beta^+_j-\beta^+_{j-1}$ $(2\le j\le N)$ and $\delta_{N+1}=-\beta^+_N$. Further, $(\omega_{N+2},\delta_{N+2})=(\theta,\delta_*)$. Finally, $\omega_j=-1$ for $j\ge N+3$, with $\delta_{N+3}=\beta^-_1$ and $\delta_{N+2+k}=\beta^-_k-\beta^-_{k-1}$ $(2\le k\le N)$.
\item Registration: $j_i:=i$ for $0\le i\le N$ and $j_i:=i+1$ for $N+1\le i\le2N+1$. Step $N+1$ is not registered.
\item Outer coefficients: $\widehat\gamma_0=\gamma$, $\widehat\gamma_i=\alpha_i$ $(1\le i\le N)$, $\widehat\gamma_{N+1}=-\gamma A_*$, and $\widehat\gamma_{N+1+k}=-\alpha_kA_*$ $(1\le k\le N)$.
\item Path maps, with $\mathcal D_{0,\eta}=\id$:
\[
\mathcal D_{j,\eta}:=F_{\omega_j,\eta}(\delta_j)\circ\dots\circ F_{\omega_1,\eta}(\delta_1),\qquad\mathcal C_{j,\eta}:=\chi\circ\mathcal D_{j,\eta}\circ\chi^{-1}\qquad(0\le j\le2N+2).
\]
\end{enumerate}
\end{cons}

\begin{lem}\label{lem:pathphase}
For every $u\in\Rr$,
\[
\sum_{i=0}^{2N+1}\widehat\gamma_i\exp\Big(\sum_{k=1}^{j_i}\delta_k e^{\omega_k u}\Big)=\psi_\theta(u).
\]
\end{lem}

\begin{proof}
Let $\beta^{\pm}_0:=0$ and $S_j(u):=\sum_{k=1}^{j}\delta_ke^{\omega_ku}$. The increments telescope within each block:
\[
S_j(u)=\beta^+_je^{u}\quad(0\le j\le N),\qquad S_{N+2+k}(u)=\delta_*e^{\theta u}+\beta^-_ke^{-u}\quad(0\le k\le N).
\]
Since $j_i=i$ for $i\le N$ and $j_{N+1+k}=N+2+k$ for $0\le k\le N$, the left-hand side equals $\Phi_+(u)-A_*e^{\delta_*e^{\theta u}}\Phi_-(u)$, which is $\psi_\theta(u)$ by \eqref{eq:Astar}.
\end{proof}

\begin{lem}\label{lem:pathfield}
For $0\le j\le 2N+2$ let $G_j:=\sum_{k=1}^{j}\delta_kw_{\omega_k}$. Then $\partial_\eta\mathcal D_{j,\eta}|_{\eta=0}=-G_j$ and $\partial_\eta\mathcal C_{j,\eta}(z)|_{\eta=0}=-\widehat G_j(z)$ with $\widehat G_j(\chi(x))=\chi'(x)G_j(x)$. In particular,
\[
\frac{\widehat G_j(\chi(x))}{\chi(x)}=\varpi_L(x)\,\frac{G_j(x)}{w(x)} .
\]
\end{lem}

\begin{proof}
By continuity in $(\eta,x)$ of all prefixes, the set of $(\eta,x)$ at which the word
\[
F_{\omega_1,\eta}(\delta_1)\ast\dots\ast F_{\omega_j,\eta}(\delta_j)
\]
is legal on $\{x\}$ is open. Since $F_{\omega,0}(\delta)=\id$ on $O_\omega$, it contains $\{0\}\times U_1$ for an open set $U_1$ with
\[
I\subset U_1\subset\mathcal D_{B\ast B^-}\cap\bigcap_\omega O_\omega,
\]
and $(\eta,x)\to\mathcal D_{j,\eta}(x)$ is real analytic on it. Induction on $j$, with $\mathcal D_{j-1,0}=\id$ on $U_1$ and \eqref{eq:movinginput} for $x(\eta)=\mathcal D_{j-1,\eta}(x)$, shows
\[
\partial_\eta\mathcal D_{j,\eta}(x)|_{\eta=0}=-G_{j-1}(x)-\delta_jw_{\omega_j}(x)=-G_j(x).
\]
As $\chi^{-1}(\chi(x))=x$ does not depend on $\eta$, the chain rule implies
\[
\partial_\eta\mathcal C_{j,\eta}(\chi(x))|_{\eta=0}=-\chi'(x)G_j(x).
\]
The last identity is $\chi'/\chi=\varpi_L/w$ (Lemma~\ref{lem:weight}(2)).
\end{proof}
\section{The compiler II: outer weights, folding, terminal closure}\label{sec:compilerB}

For compact $E\subset\Rr$ and $f$ of class $C^1$ near $E$ let $\|f\|_{C^1(E)}:=\sup_E|f|+\sup_E|f'|$.

\subsection{Outer weights and the approximation of the phase}\label{sec:outer}

\begin{lem}\label{lem:outerweight}
Let $0<z_-<z_+$, $k\in\mathbb N$, and for $i<k$ and $h\in\mathbb N^+$ let $G_i$ and $C_{i,h}$ be $C^2$ on a neighbourhood of $[z_-,z_+]$. Assume that there exist $A\ge0$ and $h_0\in\mathbb N^+$ with
\[
\big|\partial_z^{\iota}\big(C_{i,h}(z)-z+G_i(z)/h\big)\big|\le A/h^{2}\qquad(\iota\in\{0,1,2\},\ h\ge h_0,\ i<k,\ z\in[z_-,z_+]).
\]
Then there exist $h_1\in\mathbb N^+$ and $K_1\ge0$ such that for $h\ge h_1$, $i<k$, and $z\in[z_-,z_+]$: $C_{i,h}(z)>0$, $C_{i,h}'(z)>0$, and $\mathcal W_{i,h}(z):=C_{i,h}(z)^{-h}/(z^{-h}C_{i,h}'(z))$ satisfies
\[
\big|\mathcal W_{i,h}(z)-e^{G_i(z)/z}\big|+\big|\mathcal W_{i,h}'(z)-\big(e^{G_i/z}\big)'(z)\big|\le K_1/h .
\]
\end{lem}

\begin{proof}
Let $p_i:=G_i/z$ and
\begin{gather*}
c_0:=\max_{i<k}\sup_{[z_-,z_+]}(|p_i|+|p_i'|+|G_i'|+|G_i''|),\qquad c_1:=A/z_-+A/z_-^2,\\
c_2:=c_1+2(c_0+c_1)^2+2(c_0+A),\quad\text{and}\quad h_1:=\max\big(h_0,\lceil 2(c_0+c_1)\rceil,\lceil 2(c_0+A)\rceil\big).
\end{gather*}
For $h\ge h_1$ write on $[z_-,z_+]$, dropping $i$,
\[
\frac{C_{h}}{z}=1+u,\quad u:=-\frac ph+r,\quad r:=\frac{C_{h}-z+G/h}{z},\qquad C_{h}'=1+v .
\]
By hypothesis and $z\ge z_-$, $|r|,|r'|\le c_1/h^2$ and $|v|,|v'|\le(c_0+A)/h\le\frac12$, so $|u|,|u'|\le(c_0+c_1)/h\le\frac12$. Thus $C_h>0$, $C_h'>0$, and $\log\mathcal W_{h}=-h\log(1+u)-\log(1+v)$ implies
\[
\log\mathcal W_{h}-p=-hr-h\big(\log(1+u)-u\big)-\log(1+v),\qquad
(\log\mathcal W_{h})'-p'=-hr'+\frac{h\,uu'}{1+u}-\frac{v'}{1+v}.
\]
By $|\log(1+t)-t|\le2t^2$ and $|\log(1+t)|\le2|t|$ for $|t|\le\frac12$, both right-hand sides are bounded by $c_2/h$. With $c_3:=e^{c_0+c_2}$, the mean value theorem yields $|\mathcal W_h-e^p|\le c_2c_3/h$, and
\[
\mathcal W_h'-(e^p)'=\mathcal W_h\big((\log\mathcal W_h)'-p'\big)+\big(\mathcal W_h-e^p\big)p'
\]
is bounded by $c_2c_3(1+c_0)/h$. The lemma holds with $K_1:=c_2c_3(2+c_0)$.
\end{proof}

\begin{prop}\label{prop:approx}
Let $\ell,\theta,\Lambda$ and $\widehat\psi$ be as in Section~\ref{sec:phase}. For every $\epsilon>0$ there exist $m\in\mathbb N^+$, $L>0$ and $h_0\in\mathbb N^+$ such that, with $\chi:=\chi_L\circ B$ and $I':=[\chi(x_-),\chi(x_+)]$, for every integer $h\ge h_0$:
\begin{enumerate}
\item There exists an open set $U$ with $I'\subset U\subset(1,\infty)$ such that for $0\le j\le2N+2$ the words
\[
\chi^{-1},F_{\omega_1,1/h}(\delta_1),\dots,F_{\omega_j,1/h}(\delta_j),\chi
\]
are legal, each on the image of $U$ under the composition of the words before it, and likewise the words of the inverse path, and each $\mathcal C_i:=\mathcal C_{j_i,1/h}$ $(0\le i\le2N+1)$ is an analytic strictly increasing diffeomorphism of $U$ onto an open subset of $(0,\infty)$, with $\mathcal C_0=\id$.
\item The function
\[
Q_h(z):=\sum_{i=0}^{2N+1}\widehat\gamma_i\,\mathcal C_i(z)^{-h}/\big(z^{-h}\mathcal C_i'(z)\big)
\]
is $C^1$ on $U$, $\sup_{I}|Q_h\circ\chi-\widehat\psi|<\epsilon$, and $\sup_I|(Q_h\circ\chi)'-\widehat\psi'|<\epsilon$.
\item $Q_h$ has $RM$ pairwise distinct zeros in $(\chi(x_-),\chi(x_+))$ at which $Q_h'\ne0$.
\end{enumerate}
\end{prop}

\begin{proof}
By Lemmas~\ref{lem:phase} and~\ref{lem:rootmult}, $\widehat\psi$ is $C^1$ near $I$ with $RM$ simple zeros in $(x_-,x_+)$. Let $\epsilon_1$ be the tolerance of Lemma~\ref{lem:stability} for $f=\widehat\psi$ on $I$, and $\epsilon':=\min(\epsilon,\epsilon_1)/3$.

\emph{Choice of $m$.} By Lemma~\ref{lem:transport}, $w_\omega/w\to e^{\omega\Lambda V/B}$ in $C^1(I)$ as $m\to\infty$. Hence
\[
G_{j}/w=\sum_{k\le j}\delta_kw_{\omega_k}/w\to\sum_{k\le j}\delta_ke^{\omega_k\Lambda V/B}.
\]
As finite sums and $\exp$ preserve uniform $C^1$ convergence on $I$, Lemma~\ref{lem:pathphase} with $u=\Lambda V(x)/B(x)$ yields
\[
S_m:=\sum_{i=0}^{2N+1}\widehat\gamma_i\exp\big(G_{j_i}/w\big)\longrightarrow\psi_\theta(\Lambda V/B)=\widehat\psi\quad\text{in }C^1(I).
\]
Fix $m\ge m_0$ with $\|S_m-\widehat\psi\|_{C^1(I)}<\epsilon'$.

\emph{Choice of $L$.} By Lemma~\ref{lem:pathfield}, $\widehat G_j(\chi)/\chi=\varpi_L\,G_j/w$ on $I$, where $G_j/w$ does not depend on $L$ and $\|\varpi_L-1\|_{C^1(I)}\le 2C/L$ (Lemma~\ref{lem:weight}(2)). Fix $L$ with $\|S_{m,L}-S_m\|_{C^1(I)}<\epsilon'$, where
\[
S_{m,L}:=\sum_i\widehat\gamma_i\exp\big(\widehat G_{j_i}(\chi)/\chi\big).
\]

\emph{Choice of $h_0$.} The words in (1) and in the inverse path have coefficients affine in $\eta:=1/h$, so the set of $(\eta,z)$ at which all are legal is open and the compositions are jointly analytic on it (Lemma~\ref{lem:legal}(1)). It contains $\{0\}\times I'$ (at $\eta=0$ the intermediate points are in $I$, $B(I)$, or $W_{\xi_\omega}(I)$), hence, by compactness, $(-\eta_0,\eta_0)\times U$ for some $\eta_0>0$ and an open interval $U$ with $I'\subset U\subset(1,\infty)$. On $U$, $\mathcal C_{j,0}=\id$, and $\partial_\eta\mathcal C_{j,\eta}|_{\eta=0}=-\widehat G_j$ (Lemma~\ref{lem:pathfield}). Taylor's theorem in $\eta$ \cite{Rud76} for $\partial_z^\iota\mathcal C_{j,\eta}(z)$ yields $A\ge0$ with
\[
\big|\partial_z^{\iota}\big(\mathcal C_{j,\eta}(z)-z+\eta\,\widehat G_j(z)\big)\big|\le A\eta^2\qquad(\iota\le2,\ |\eta|\le\eta_0/2,\ j\le2N+2,\ z\in I').
\]
Lemma~\ref{lem:outerweight} with $[z_-,z_+]=I'$, $G_i=\widehat G_{j_i}$, and $C_{i,h}=\mathcal C_{j_i,1/h}$ for $h\ge2/\eta_0$ ($C_{i,h}:=\id$ otherwise, and $i\le2N+1$) and $h_0$ replaced by $\lceil2/\eta_0\rceil$ yields $h_1$ and $K_1$. Since $Q_h=\sum_i\widehat\gamma_i\mathcal W_{i,h}$, the chain rule yields
\[
\|Q_h\circ\chi-S_{m,L}\|_{C^1(I)}\le\Big(1+\sup_I\chi'\Big)\Big(\sum_{i=0}^{2N+1}|\widehat\gamma_i|\Big)\frac{K_1}{h}\qquad(h\ge h_1).
\]
Let $h_0\ge\max(h_1,2/\eta_0)$ be an integer for which the right-hand side is $<\epsilon'$.

Let $h\ge h_0$. (1) On $U$, $\mathcal C_i$ is analytic with $\mathcal C_i'>0$ (Lemma~\ref{lem:legal}(1)) and inverse given by the inverse path, and $\mathcal C_0=\chi\circ\chi^{-1}=\id$ on $U$ by Lemma~\ref{lem:legal}(2) with $O=\chi((x_-,x_+))$. (2) $Q_h$ is $C^1$ on $U$ as $\mathcal C_i,\mathcal C_i',z>0$, and $\|Q_h\circ\chi-\widehat\psi\|_{C^1(I)}<3\epsilon'\le\epsilon$. (3) As $3\epsilon'\le\epsilon_1$, Lemma~\ref{lem:stability} yields $RM$ simple zeros $q$ of $Q_h\circ\chi$ in $(x_-,x_+)$. As $\chi'>0$ on $I$, the points $\chi(q)$ are $RM$ distinct zeros of $Q_h$ in $(\chi(x_-),\chi(x_+))$ with $Q_h'(\chi(q))=(Q_h\circ\chi)'(q)/\chi'(q)\neq0$.
\end{proof}

\subsection{The interleaved family}\label{sec:interleaved}

\begin{cons}\label{cons:interleaved}
Let $h$ and $\mathcal C_0,\dots,\mathcal C_{2N+1}$ be as in Proposition~\ref{prop:approx}. Let $\Theta_{i,u}:=\Psi^{(h)}_{u,\widehat\gamma_i}$ $(0\le i\le 2N+1)$, $K:=\mathcal C_{2N+1}$ and
\[
\mathcal H_u:=\Theta_{2N+1,u}\circ\big(\mathcal C_{2N+1}\circ\mathcal C_{2N}^{-1}\big)\circ\Theta_{2N,u}\circ\dots\circ\Theta_{1,u}\circ\big(\mathcal C_{1}\circ\mathcal C_{0}^{-1}\big)\circ\Theta_{0,u}.
\]
\end{cons}

\begin{lem}\label{lem:chainvar}
Let $h\in\mathbb N^+$, $P\in\mathbb N$, $\widehat\gamma_0,\dots,\widehat\gamma_P\in\Rr$, $U\subset(0,\infty)$ open, and let $\mathcal C_0=\id,\mathcal C_1,\dots,\mathcal C_P$ be $C^2$ diffeomorphisms of $U$ onto open subsets of $(0,\infty)$. Define $\mathcal H_u$, $K:=\mathcal C_P$, and $Q_h$ by Construction~\ref{cons:interleaved} and Proposition~\ref{prop:approx}(2) with $P$ in place of $2N+1$. Then:
\begin{enumerate}
\item For every compact $J\subset U$ there exist $\epsilon>0$ and an open set $J\subset U_J\subset U$ such that $\mathcal H_u$ is defined on $U_J$ for $|u|<\epsilon$. Moreover, $\mathcal H_0=K$ on $U$.
\item $\partial_u\mathcal H_u(z)|_{u=0}=K'(z)\,z^{-h}\,Q_h(z)$ for $z\in U$.
\item For $z\in U$: the function $z\to\partial_u\mathcal H_u(z)|_{u=0}$ vanishes at $z$ with nonzero derivative if and only if $Q_h(z)=0\ne Q_h'(z)$.
\end{enumerate}
\end{lem}

\begin{proof}
Let $Z_0:=\Theta_{0,u}$ and
\[
Z_{i+1}:=\Theta_{i+1,u}\circ\mathcal C_{i+1}\circ\mathcal C_i^{-1}\circ Z_i,
\]
so $\mathcal H_u=Z_P$. Since $\Theta_{i,0}=\id$ (Lemma~\ref{lem:powerflow}) and $\mathcal C_0=\id$, induction yields $Z_i|_{u=0}=\mathcal C_i$ on $U$. In particular $\mathcal H_0=K$. At $u=0$ the intermediate points are in open domains and the radicands are $\mathcal C_i(z)^{h+1}>0$. Continuity and compactness of $\{0\}\times J$ yield (1). Let $v_i(z):=\partial_uZ_i(z)|_{u=0}$. Lemma~\ref{lem:powerflow} (with $h,u,\widehat\gamma_i$ in place of $m,\eta,\delta$) and
\[
(\mathcal C_{i+1}\circ\mathcal C_i^{-1})'(\mathcal C_i(z))=\mathcal C_{i+1}'(z)/\mathcal C_i'(z)
\]
yield
\[
v_0=\widehat\gamma_0z^{-h},\qquad v_{i+1}=\widehat\gamma_{i+1}\,\mathcal C_{i+1}^{-h}+\frac{\mathcal C_{i+1}'}{\mathcal C_i'}\,v_i,
\qquad\text{hence}\qquad v_i=\mathcal C_i'\sum_{k=0}^{i}\widehat\gamma_k\frac{\mathcal C_k^{-h}}{\mathcal C_k'} .
\]
For $i=P$ this is (2). The factor $z\to K'(z)z^{-h}$ is $C^1$ and nonzero on $U$, and $Q_h$ is $C^1$ on $U$ since the $\mathcal C_k$ are $C^2$. At a zero of $Q_h$ the product rule yields $(K'z^{-h}Q_h)'=K'z^{-h}Q_h'$. This proves (3).
\end{proof}

\subsection{Budgets}\label{sec:budget}

In the words below, $\chi$ and $\chi^{-1}$ denote the words $B\ast\chi_L$ and $\chi_L^{-1}\ast B^{-}$ of lengths at most $b+2$ (Lemma~\ref{lem:weight}(1)), and $\Theta_{i,u}$ the two-segment word of Lemma~\ref{lem:powerflow}.

\begin{lem}\label{lem:budget}
Let $m,L$ be as in Proposition~\ref{prop:approx}. There exists $h_1\in\mathbb N^+$ such that for every integer $h\ge h_1$ there exist $\epsilon_2>0$ and an open set $U_2\supset I'$ with the following properties.
\begin{enumerate}
\item For $|u|\le\epsilon_2$, $\mathcal H_u$ coincides on $U_2$ with a word $\mathcal W_u$, where $u\to\mathcal W_u$ is a word family with analytic coefficients, $|\mathcal W_u|\le(4N+4)n+(4N+2)b+16N+12$ and $\Pi(\mathcal W_u)=1$. A word family $\mathcal W^{-}_u$ with $|\mathcal W^{-}_u|\le(4N+4)n+(4N+2)b+16N+12$ is legal on $\mathcal W_u(U_2)$ and satisfies $\mathcal W^{-}_u\circ\mathcal W_u=\id$ on $U_2$.
\item On $U_2$,
\[
K=\chi\circ F_{-1,1/h}(\beta^-_N)\circ F_{\theta,1/h}(\delta_*)\circ\chi^{-1},
\]
and
\[
\chi\circ F_{-1,1/h}(\beta^-_N)\circ F_{\theta,1/h}(\delta_*)\circ\chi^{-1}\qquad\text{and}\qquad\chi\circ F_{\theta,1/h}(-\delta_*)\circ F_{-1,1/h}(-\beta^-_N)\circ\chi^{-1}
\]
are words of length at most $4n+2b+8$ with exponent product $1$, legal on $U_2$ and on $K(U_2)$ respectively, and the second composed with the first is $\id$ on $U_2$.
\end{enumerate}
\end{lem}

\begin{proof}
For $0\le i\le2N$ let
\[
T_i:=\chi^{-1}\ast F_{\omega_{j_i+1},1/h}(\delta_{j_i+1})\ast\dots\ast F_{\omega_{j_{i+1}},1/h}(\delta_{j_{i+1}})\ast\chi.
\]
As $\mathcal D_{j_{i+1},\eta}$ extends $\mathcal D_{j_i,\eta}$ by the steps $j_i+1,\dots,j_{i+1}$, the inner $\chi^{-1}\circ\chi$ cancel and $T_i=\mathcal C_{i+1}\circ\mathcal C_i^{-1}$ on $\mathcal C_i(U)$. Let
\[
\mathcal W_u:=\Theta_{0,u}\ast T_0\ast\Theta_{1,u}\ast T_1\ast\dots\ast T_{2N}\ast\Theta_{2N+1,u},
\]
and let $\mathcal W^-_u$ be the reversed list with $\Theta_{i,-u}$ in place of $\Theta_{i,u}$ and each $T_i$ replaced by
\[
\chi^{-1}\ast(\text{block }i\text{ reversed, increments negated})\ast\chi.
\]
The blocks contain $2N+2$ steps in total ($j_{N+1}-j_N=2$, all other gaps $1$), so
\[
|\mathcal W_u|\le(2N+2)(2n+2)+(2N+1)(2b+4)+2(2N+2)=(4N+4)n+(4N+2)b+16N+12,
\]
and likewise for $\mathcal W^-_u$. The coefficients are affine in $u$ (only $(h+1)\widehat\gamma_iu$ in $\Theta_{i,u}$ varies). Each $\Theta_{i,u}$ and each $F_{\omega,\eta}(\delta)$ has exponent product $1$, and $\Pi(\chi)\Pi(\chi^{-1})=\Pi(B)\Pi(B^-)=1$ by Lemmas~\ref{lem:weight}(1) and~\ref{lem:legal}(2). Thus $\Pi(\mathcal W_u)=1$.

For (2), the group law of $F_{\omega,\eta}$ in $\delta$ within each block, with $\sum_{k=1}^{N+1}\delta_k=0$ and $\sum_{k=1}^{N}\delta_{N+2+k}=\beta^-_N$, and $F_{1,\eta}(0)=W^-_{\xi_1}\circ W_{\xi_1}=\id$ near $I$ (Lemma~\ref{lem:legal}(2)), yield
\[
\mathcal D_{2N+2,\eta}=F_{-1,\eta}(\beta^-_N)\circ F_{\theta,\eta}(\delta_*)
\]
near $I$ for small $|\eta|$. As $K=\mathcal C_{j_{2N+1},1/h}=\mathcal C_{2N+2,1/h}$, this proves the identity in (2). The words in (2) have length at most $2(2n+2)+2(b+2)=4n+2b+8$ and exponent product $1$, and by the group law the second composed with the first is $\id$.

As for $h_0$ in Proposition~\ref{prop:approx}, compactness yields $h_1\ge h_0$ such that for $h\ge h_1$, after shrinking $U$ around $I'$, the folding identity holds on $U$, $T_i$ is legal on $\mathcal C_i(U)$, the reversal of $T_i$ on $\mathcal C_{i+1}(U)$, and the words of (2) on $U$ and $K(U)$, respectively. Let $\epsilon,U_J$ be as in Lemma~\ref{lem:chainvar}(1) with $J=I'$ (so $\mathcal C_i^{-1}$ acts only on $\mathcal C_i(U)$), $\epsilon_2:=\epsilon/2$ and $U_2:=U_J$. Evaluating atom by atom along $\mathcal H_u(z)$ ($z\in U_2$), the group laws of Lemma~\ref{lem:powerflow} and of $F_{\omega,1/h}$ yield legality, $\mathcal W_u=\mathcal H_u$ and $\mathcal W^-_u\circ\mathcal W_u=\id$ on $U_2$.
\end{proof}

\subsection{Proof of Theorem~\ref{thm:compiler}}\label{sec:compilerproof}

If $RM=0$, $W$ with the empty list of points is an element of $\mathcal F(n,b,0)$. Take $n'=n$, $b'=b$. Otherwise $M\ge1$, so $x_-<x_+$. The choices are made in the following order.
\begin{enumerate}
\item $\ell,\theta$ (Lemma~\ref{lem:phase}).
\item $\Lambda$ (Lemma~\ref{lem:rootmult} with $\psi=\psi_\theta$).
\item $m,L,h_0$ (Proposition~\ref{prop:approx} with $\epsilon=1$).
\item $h_1$ (Lemma~\ref{lem:budget}) and $h:=\max(h_0,h_1)$.
\item $\epsilon_2,U_2$ (Lemma~\ref{lem:budget}).
\end{enumerate}

With $n',b'$ the bounds of Theorem~\ref{thm:compiler}, the following data form an element of $\mathcal F(n',b',RM)$.
\begin{enumerate}
\item The interval $I'$, with $1<\chi(x_-)<\chi(x_+)$ by Lemma~\ref{lem:weight}, and the radius $\epsilon_2$.
\item The families $\mathcal W_u,\mathcal W^-_u$ of Lemma~\ref{lem:budget}(1).
\item The words of Lemma~\ref{lem:budget}(2) as base words. The first of them is $K=\mathcal H_0=\mathcal W_0$ on $U_2\supset I'$ (Lemma~\ref{lem:chainvar}(1)).
\item The $RM$ zeros of $Q_h$ of Proposition~\ref{prop:approx}(3), in increasing order.
\end{enumerate}
As $\partial_u\mathcal W_u|_{u=0}=\partial_u\mathcal H_u|_{u=0}$ on the open set $U_2$, Lemma~\ref{lem:chainvar}(3), with hypotheses from Proposition~\ref{prop:approx}(1), shows that the points of (4) are simple zeros of the first variation of $\mathcal W_u$. \qed

\subsection{Terminal closure}\label{sec:terminal}

\begin{proof}[Proof of Proposition~\ref{prop:terminal}]
As $M\ge1$, $x_-<x_+$. By Lemma~\ref{lem:legal}(1), $B'>0$ near $I$, so $0<B(x_-)<B(x_+)$ and $J_B:=[B(x_-),B(x_+)]=B(I)$, and $B^-$ is legal on $J_B$ with $B^-(J_B)=I$ and $(B^-)'>0$ near $J_B$. For $|u|\le\xi_0$, $P_u:=B^-\ast W_u$ is legal on $J_B$, $|P_u|\le n+b$, and $\Pi(P_u)=\Pi(B^-)\Pi(W_u)=1$ (Lemma~\ref{lem:legal}(2) with $O=(x_-,x_+)$). For $y\in J_B$, $P_0(y)=B(B^-(y))=y$ because $W_0=B$ on $I$. By Lemma~\ref{lem:legal}(2) with $O=(B(x_-),B(x_+))$, $P_0=\id$ on an open set containing $J_B$.

As the coefficients are analytic, Lemma~\ref{lem:legal}(1) and compactness yield $0<r\le\xi_0$ such that $(u,y)\to P_u(y)$ is $C^2$ near $[-r,r]\times J_B$. Let $D_u:=(P_u-\id)/u$ and $f:=V\circ B^-$. Then $\partial_uP_u(y)|_{u=0}=f(y)$. As $P_0=\id$ near $J_B$, the mean value theorem \cite{Rud76} in $u$ yields $D_u(y)=\partial_uP_{tu}(y)$ and $D_u'(y)=\partial_u\partial_yP_{su}(y)$ for some $t,s\in(0,1)$. By uniform continuity on $[-r,r]\times J_B$, $\|D_u-f\|_{C^1(J_B)}\to0$ as $u\to0$.

The points $B(x_1)<\dots<B(x_M)$ of $(B(x_-),B(x_+))$ satisfy
\[
f(B(x_k))=V(x_k)=0,\qquad f'(B(x_k))=V'(x_k)\,(B^-)'(B(x_k))\neq0 .
\]
With the tolerance $\epsilon_3$ of Lemma~\ref{lem:stability} for $f$ on $J_B$ at $B(x_1),\dots,B(x_M)$, let $\rho\in(0,r]$ satisfy $\|D_u-f\|_{C^1(J_B)}<\epsilon_3$ for $0<u<\rho$, and let $u:=\rho/2$. Lemma~\ref{lem:stability} yields pairwise distinct $p_1,\dots,p_M\in(B(x_-),B(x_+))$ with $D_u(p_k)=0\neq D_u'(p_k)$. Since $u\neq0$, $P_u(p_k)=p_k$ and $P_u'(p_k)=1+uD_u'(p_k)\neq1$, with $P:=P_u$.
\end{proof}
\section{The host I: polycycle, field, perturbations, full rank}\label{sec:hostA}

\subsection{The Chebyshev polycycle}\label{sec:polycycle}

\begin{setup}\label{setup:host}
Let $s\ge2$ be an integer and $S:=4s^2+2$. Let $T_s$ be the Chebyshev polynomial of degree $s$, $f:=T_s/2$, $y_0:=1/4$, $y_j:=9/10+j/(40s+20)$ $(1\le j\le 2s)$, $\mathcal R(y):=\prod_{j=1}^{2s}(y-y_j)$,
\[
M_{\mathcal R}:=\sup_{[-2,2]}(|\mathcal R|+|\mathcal R'|+|\mathcal R''|),
\]
$\varepsilon_0:=10^{-6}/(1+M_{\mathcal R})$, $\phi(y):=y^2+\varepsilon_0\mathcal R(y)$, and
\[
A:=f(x)^2+y^2-1,\qquad B:=f(x)^2+\phi(y)-1,\qquad C:=y_0-y,
\]
$\Omega:=\{A<0,\ B<0,\ C<0\}$, $\Gamma:=\partial\Omega$, $J(u,v):=(-v,u)$.
\end{setup}

\begin{lem}\label{lem:height}
In Set-up~\ref{setup:host}:
\begin{enumerate}
\item We have
\[
M_{\mathcal R}=\mathcal R(-2)-\mathcal R'(-2)+\mathcal R''(-2),
\]
and $\varepsilon_0|\mathcal R^{(\iota)}|<10^{-6}$ on $[-2,2]$ for $\iota=0,1,2$.
\item $9/10<y_1<\dots<y_{2s}<19/20$, $\mathcal R(y_0)>0$, and $\mathcal R'(y_j)\ne0$.
\item $\phi''>1$ on $\Rr$ and $\phi'>0$ on $[y_0,\infty)$.
\item $y_0^2<\phi(y_0)<1/4$ and $3/4<\phi(89/100)<\phi(97/100)<1<\phi(1)$.
\end{enumerate}
\end{lem}

\begin{proof}
For $|y|\le2$, $|y-y_j|\le|-2-y_j|$, and all factors of $\mathcal R$ are negative at $-2$. Termwise in the product rule,
\[
|\mathcal R^{(\iota)}(y)|\le(-1)^\iota\mathcal R^{(\iota)}(-2),
\]
hence (1). Item (2) uses $\mathcal R'(y_j)=\prod_{i\ne j}(y_j-y_i)$. Item (3) uses (1), $\mathcal R''>0$ on $\{|y|>2\}$ and $\phi'(y_0)=1/2+\varepsilon_0\mathcal R'(y_0)$. Item (4) uses (1), (3), $\mathcal R(y_0),\mathcal R(1)>0$ and $1/16+10^{-6}<1/4$, $0.7921-10^{-6}>3/4$, $0.9409+10^{-6}<1$.
\end{proof}

\begin{prop}\label{prop:polycycle}
In Set-up~\ref{setup:host}, $\Omega$ is a bounded Jordan domain and $\Gamma$ contains exactly $S$ points at which two of $A,B,C$ vanish (the corners), namely the $4s^2$ points with $A=B=0$, which are on the lines $y=y_j$ $(1\le j\le2s)$, $2s$ on each, with $|x|<1$, and the points $(\pm b_*,y_0)$, $b_*>1$, with $B=C=0$. At each corner the vanishing factors have linearly independent gradients and the third factor is negative. The corners can be numbered $p_0=(-b_*,y_0)$, $p_1=(b_*,y_0)$, $p_2,\dots,p_{S-1}$ so that $\Gamma$ minus the corners is the disjoint union of open regular smooth arcs $\Gamma_k$ from $p_k$ to $p_{k+1}$ $(k\in\ZZ/S)$, and on $\Gamma_k$ exactly one factor $F_k\in\{A,B,C\}$ vanishes, with $F_0=C$, $F_k=B$ for odd $k$ and $F_k=A$ for even $k\ge2$. In particular $\#\{k:F_k=A\}=2s^2$, $\#\{k:F_k=B\}=2s^2+1$ and $\#\{k:F_k=C\}=1$.
\end{prop}

\begin{proof}
By $T_s(\cos\theta)=\cos s\theta$ and $T_s(\cosh\theta)=\cosh s\theta$ \cite{Riv90}, for $0<u<1/4$ the equation $f(x)^2=u$ has exactly $2s$ solutions, all simple and in $(-1,1)$. For $u>1/4$ exactly two, $\pm x_u$, simple, $x_u>1$, and $\{f^2<u\}=(-x_u,x_u)$. As $B-A=\varepsilon_0\mathcal R$, Lemma~\ref{lem:height} implies the following.
\begin{enumerate}
\item $A=B=0$ iff $y=y_j$ and $f(x)^2=1-y_j^2\in(0,1/4)$.
\item $B=C=0$ iff $y=y_0$ and $f(x)^2=1-\phi(y_0)>3/4$, i.e.\ $x=\pm b_*$.
\item $A=C=0$ iff $y=y_0$ and $f(x)^2=1-y_0^2>3/4$, i.e.\ $x=\pm a_*$ for some $a_*>1$, and $B(\pm a_*,y_0)=\phi(y_0)-y_0^2>0$.
\end{enumerate} At the corners
\[
\det(\nabla A,\nabla B)=\varepsilon_0(f^2)'(x)\,\mathcal R'(y_j)\ne0,\qquad \det(\nabla B,\nabla C)=-(f^2)'(\pm b_*)\ne0,
\]
$C=y_0-y_j<0$ and $A=y_0^2-\phi(y_0)<0$. With $\phi_+^{-1}:=(\phi|_{[y_0,1]})^{-1}$, smooth near $[\phi(y_0),\phi(1)]$, monotonicity yields
\begin{equation}\label{eq:roof}
\Omega=\{|x|<b_*,\ y_0<y<r(x)\},\qquad r:=\min\big(\sqrt{1-f^2},\ \phi_+^{-1}(1-f^2)\big),
\end{equation}
with $r(\pm b_*)=y_0<r$ on $(-b_*,b_*)$. Via $(x,t)\to(x,y_0+t(r(x)-y_0))$, $\Omega$ is a bounded Jordan domain with $\Gamma=([-b_*,b_*]\times\{y_0\})\cup\operatorname{graph}r$. The branches of $r$ meet exactly at the points $A=B=0$, so between consecutive corners $\operatorname{graph}r$ is one branch, smooth up to the ends. Number the upper corners by decreasing abscissa. As $\partial_yA,\partial_yB>0$ on $\Gamma$ and $x$ decreases along $\Gamma_k$ ($k\ne0$), a tangent $\dot z\perp\nabla F_k$ equals $-(\dot x/\partial_yF_k)J\nabla F_k$. With $J\nabla C=(1,0)$,
\begin{equation}\label{eq:orient}
\text{the velocity of }\overline{\Gamma_k}\text{ from }p_k\text{ to }p_{k+1}\text{ is a positive multiple of }J\nabla F_k\quad(k\in\ZZ/S).
\end{equation}
Near $p_k$ the vanishing factors $(F,G)$ form a chart in which $\overline\Omega=\{F\le0,\,G\le0\}$, so $F_{k-1}\ne F_k$. With $F_k\in\{A,B\}$ for $k\ne0$, $F_1=B$, and $S-1$ odd, the pattern follows.
\end{proof}

\subsection{The field with given saddle ratios}\label{sec:field}

\begin{cons}\label{cons:field}
In Set-up~\ref{setup:host} let $\nu_0,\dots,\nu_{S-1}>0$, $|\tau|<\nu_{S-1}/2$, and $\nu_k(\tau):=\nu_k+\tau$ if $k=S-1$, $\nu_k(\tau):=\nu_k$ otherwise. For $F\in\{A,B,C\}$ let $w_F(\tau)$ be the polynomial
\[
w_F(\tau)(x,y):=\sum_{j=0}^{2s}L_j(y)^2\sum_{\xi\in N_j}d_F(\tau)(\xi,y_j)\,L_{j,\xi}(x)^2,
\]
where $L_j$ is the Lagrange basis of the heights $y_0,\dots,y_{2s}$, $N_0$ is the set of the four real zeros of $A(\cdot,y_0)B(\cdot,y_0)$, $N_j$ $(j\ge1)$ the $2s$ real zeros of $A(\cdot,y_j)$, $L_{j,\xi}$ the Lagrange basis of $N_j$, and $d_F(\tau)(p):=\nu_k(\tau)$ if $p=p_k$ and $F=F_{k-1}$, $d_F(\tau)(p):=1$ otherwise. Define
\[
X_0(\tau):=w_A(\tau)\,BC\,J\nabla A+w_B(\tau)\,AC\,J\nabla B+w_C(\tau)\,AB\,J\nabla C .
\]
\end{cons}

\begin{lem}\label{lem:interp}
In Construction~\ref{cons:field}, $w_F(\tau)>0$ on $\Rr^2$, $\deg w_F(\tau)\le8s-2$, $w_{F_{k-1}}(\tau)(p_k)=\nu_k(\tau)$, and $w_{F_k}(\tau)(p_k)=1$ ($F\in\{A,B,C\}$, $k\in\ZZ/S$).
\end{lem}

\begin{proof}
The corners are nodes ($N_0=\{\pm a_*,\pm b_*\}$), $\#N_j\le2s$, and all data are positive as $\nu_{S-1}+\tau>0$. As $\sum_jL_j=\sum_\xi L_{j,\xi}=1$, some $L_j(y)$ and some $L_{j,\xi}(x)$ are nonzero, so $w_F(\tau)>0$. The Kronecker property and $F_{k-1}\ne F_k$ yield the values, and $\deg w_F(\tau)\le4s+2(2s-1)$.
\end{proof}

\begin{cons}\label{cons:pert}
Let $y^A_l:=(y_{2l-1}+y_{2l})/2$ $(1\le l\le s)$, $y^B_0:=89/100$, $y^B_s:=97/100$, $y^B_l:=(y_{2l}+y_{2l+1})/2$ $(1\le l\le s-1)$, $U_A(y):=\prod_{l=1}^s(y-y^A_l)$, $U_B(y):=\prod_{l=0}^s(y-y^B_l)$ and
\[
Z:=(A^2+U_A^2)(B^2+U_B^2)(C^2+x^2).
\]
For $k\in\ZZ/S$ let $P_0:=1$ and, for $k\ne0$, let $P_k$ be the polynomial $\sum_l\Lambda_l(y)I_{k,l}(x)$, where $\Lambda_l$ is the Lagrange basis of the heights $y^{F_k}_l$ and $I_{k,l}$ interpolates on the real zeros of $F_k(\cdot,y^{F_k}_l)$ the indicator of $\Gamma_k$. With $K_A:=BC$, $K_B:=AC$, $K_C:=AB$, $Y_k:=K_{F_k}^2P_k\nabla F_k$, and for $T\in\Rr$, $\zeta\in\Rr^S$,
\[
X_{T,\zeta}(\tau):=(1+TZ)\,X_0(\tau)+\sum_{k\in\ZZ/S}\zeta_kY_k .
\]
\end{cons}

\begin{lem}\label{lem:samples}
\begin{enumerate}
\item $Z\ge0$, $\deg Z=8s+2$, and the zeros of $Z$ in $\Gamma$ are: exactly one point on each arc in $\{A=0\}$, one or two points on each arc in $\{B=0\}$, and $(0,y_0)$ on $\Gamma_0$, and $Z>0$ at the corners. If $e\in\Gamma_k$ and $Z(e)=0$, then $e$ is a nondegenerate quadratic zero of $Z\circ c$ for every regular $C^2$ parametrization $c$ of $\Gamma_k$.
\item $P_k=1$ at the zeros of $Z$ on $\Gamma_k$ and $P_k=0$ at the zeros of $Z$ on the arcs $\Gamma_{k'}\ne\Gamma_k$ with $F_{k'}=F_k$. Moreover, $\deg P_k\le3s-2$ if $F_k=A$ and $\deg P_k\le3s-1$ if $F_k=B$.
\end{enumerate}
\end{lem}

\begin{proof}
(1) $Z$ is a product of sums of squares, $\deg Z\le8s+2$ and $\deg Z(x,0)=8s+2$. By Lemma~\ref{lem:height}, $y_0<y^B_0<y_1<y^A_1<y_2<y^B_1<\dots<y_{2s}<y^B_s$, and $\mathcal R<0$ exactly on the bands $(y_{2l-1},y_{2l})$. As $B-A=\varepsilon_0\mathcal R$ and $A,B\le0$ on $\Gamma$, a non-corner point of $\Gamma$ at height $y>y_0$ is on an $A$-arc iff $\mathcal R(y)<0$. So the zeros of $Z$ in $\Gamma$ are $(0,y_0)$ and the crossings of $A$-arcs with $y=y^A_l$ and of $B$-arcs with $y=y^B_l$, none a corner. By $f(\cos\theta)^2=\cos^2(s\theta)/4$, $r$ of \eqref{eq:roof} is strictly monotone between consecutive extrema: $s$ maxima, of height $\phi_+^{-1}(1)>97/100$, and $s-1$ minima, of height $\phi_+^{-1}(3/4)<89/100$. Hence an $A$-arc joins heights $y_{2l-1},y_{2l}$ monotonically (one crossing of $y^A_l$), and a $B$-arc joins $y_{2l},y_{2l+1}$ or $y_0,y_1$ (one crossing), or contains an extremum (two crossings of $y^B_s$ or $y^B_0$). Along a regular $c$ on an $A$-arc, $Z\circ c=U_A(y\circ c)^2h$ with $h>0$ at the zero, and $(y\circ c)'\ne0$ there since $\nabla A\cdot c'=0$ and $\partial_xA=(f^2)'\ne0$ (as $0<1-(y^A_l)^2<1/4$). The $B$-arcs and $\Gamma_0$ ($Z\circ c=x(c)^2h$) are analogous.

(2) The zeros of $Z$ on the arcs $\Gamma_{k'}$ with $F_{k'}=F_k$ are nodes $(\xi,y^{F_k}_l)$, where $P_k$ is the indicator of $\Gamma_k$, and $\deg P_k\le(s-1)+(2s-1)$, respectively $s+(2s-1)$.
\end{proof}

\begin{prop}\label{prop:saddles}
Let $T\ge0$, $|\tau|<\nu_{S-1}/2$, and $\zeta\in\Rr^S$. Then:
\begin{enumerate}
\item $\deg X_{T,\zeta}(\tau)\le20s$.
\item Each $p_k$ is a hyperbolic saddle of $X_{T,\zeta}(\tau)$ with eigenvalues $-\mu_k<0<\lambda_k$, $\mu_k/\lambda_k=\nu_k(\tau)$, stable direction tangent to $\Gamma_{k-1}$ and unstable direction tangent to $\Gamma_k$, and $DX_{T,\zeta}(\tau)(p_k)=(1+TZ(p_k))DX_0(\tau)(p_k)$.
\item If $\zeta_k=0$ for every $k$ with $F_k=F$, then $X_{T,\zeta}(\tau)=X_{T,0}(\tau)$ on $\{F=0\}$.
\item Each $\Gamma_k$ is one orbit of $X_{T,0}(\tau)$, with $\alpha$-limit $p_k$ and $\omega$-limit $p_{k+1}$.
\end{enumerate}
\end{prop}

\begin{proof}
(1) By Lemmas~\ref{lem:interp} and~\ref{lem:samples}, $\deg X_0(\tau)\le12s-2$, $\deg Z\le8s+2$, and $\deg Y_k\le9s$.

(2) Let $G$ be the factor with $G(p_k)<0$ and $D:=\det(\nabla F_k,\nabla F_{k-1})(p_k)$. As $F_k<0$ on $\Gamma_{k-1}$ and $F_k(p_k)=0$, \eqref{eq:orient} implies $-D=\nabla F_k\cdot J\nabla F_{k-1}(p_k)\ge0$, so $GD>0$. As two factors vanish at $p_k$, $X_0(\tau)(p_k)=0$,
\[
DX_0(\tau)(p_k)z=G\,w_{F_{k-1}}(\tau)\,(\nabla F_k\cdot z)\,J\nabla F_{k-1}+G\,w_{F_k}(\tau)\,(\nabla F_{k-1}\cdot z)\,J\nabla F_k\quad\text{at }p_k ,
\]
which in coordinates $(\nabla F_k\cdot z,\nabla F_{k-1}\cdot z)$ is $\operatorname{diag}(-GD\,\nu_k(\tau),GD)$ by Lemma~\ref{lem:interp}. The axes are tangent to $\Gamma_{k-1}$, $\Gamma_k$. As $K_{F_b}(p_k)=0$, each $Y_b$ has zero first jet at $p_k$. With $X_0(\tau)(p_k)=0$ this implies the formula for $DX_{T,\zeta}(\tau)(p_k)$, and $1+TZ(p_k)\ge1$.

(3) If $F_b\ne F$, then $F$ divides $K_{F_b}$, so $Y_b=0$ on $\{F=0\}$.

(4) On $\Gamma_k$, $X_{T,0}(\tau)=(1+TZ)\,w_{F_k}(\tau)K_{F_k}J\nabla F_k$ with positive scalar factors. By \eqref{eq:orient}, $X_{T,0}(\tau)(c(r))=a(r)c'(r)$ for a regular $c:[0,1]\to\overline{\Gamma_k}$ from $p_k$ to $p_{k+1}$, with $a\in C^1$, $a>0$ on $(0,1)$, $a(0)=a(1)=0$. So $a(r)\le M\min(r,1-r)$, $r\to\int_{1/2}^ra^{-1}$ maps $(0,1)$ onto $\Rr$, and $c$ composed with the inverse is the orbit.
\end{proof}

\subsection{Separatrix gaps and full rank}\label{sec:fullrank}

\begin{lem}\label{lem:gapreg}
Let $\nu_0,\dots,\nu_{S-1}>0$. There exist $C^\infty$ embeddings $\Sigma_k:(-r_k,r_k)\to\Rr^2$ $(k\in\ZZ/S)$ with $\Sigma_k(0)\in\Gamma_k$, $\Sigma_k(t)\in\Omega$ for $t\in(0,r_k)$, independent of $T$, such that for every $T>0$ there exist an open set $\mathcal P\ni0$ in $\Rr^S\times\Rr$ with $|\tau|<\nu_{S-1}/2$ on $\mathcal P$, on which every $\Sigma_k$ is transverse to $X_{T,\zeta}(\tau)$, and $C^\infty$ functions $u_k,v_k$ on $\mathcal P$ with the following properties: for $p=(\zeta,\tau)\in\mathcal P$ the unstable separatrix of $p_k$ and the stable separatrix of $p_{k+1}$ of $X_{T,\zeta}(\tau)$ that continue $\Gamma_k$ meet $\Sigma_k$ first at $\Sigma_k(u_k(p))$ and $\Sigma_k(v_k(p))$. The function $\Delta_k:=u_k-v_k$ vanishes on $\{\zeta=0\}\cap\mathcal P$. Moreover, $\partial_{\zeta_b}\Delta_k(0)=0$ if $F_b\ne F_k$.
\end{lem}

\begin{proof}
Let $m_k\in\Gamma_k$ and $\Sigma_k(t):=m_k-t\nabla F_k(m_k)$. As $(F_k\circ\Sigma_k)'(0)<0$, for small $r_k$, $\Sigma_k(t)\in\Omega$ for $t>0$ and $\Sigma_k((-r_k,r_k))\cap\Gamma=\{m_k\}$. Fix $T>0$. By Proposition~\ref{prop:saddles}(4), $X_{T,0}(0)(m_k)\ne0$ is tangent to $\Gamma_k$, hence transverse to $\Sigma_k$. Invariant manifolds with parameters at the saddles of Proposition~\ref{prop:saddles}(2) \cite[Ch.~IX]{Har02}, \cite[Ch.~2]{PdM82} and smooth dependence of the flow on initial data \cite[Ch.~V]{Har02} yield $\mathcal P$ and $u_k,v_k$ with $u_k(0)=v_k(0)=0$. If $X_{T,\zeta}(\tau)=X_{T,0}(\tau)$ on $\Gamma_k$, then $\Gamma_k$, an orbit from $p_k$ to $p_{k+1}$, is both continued separatrices (local separatrices are unique) and hits $\Sigma_k$ only at $m_k$. So $\Delta_k(\zeta,\tau)=0$. By Proposition~\ref{prop:saddles}(3) with $F=F_k$, this holds for $\zeta=0$ and for $(\zeta,\tau)=(te_b,0)$ with $F_b\ne F_k$, whence $\partial_{\zeta_b}\Delta_k(0)=0$.
\end{proof}

\begin{lem}\label{lem:concentration}
Let $a_-<a_+$, $\mu$ a finite measure on $[a_-,a_+]$, $\mathcal Z$ continuous and $\ge0$ on $[a_-,a_+]$ whose zero set is a finite nonempty set $E\subset(a_-,a_+)$, and assume that near each $e\in E$, $\mathcal Z$ is $C^2$ with $\mathcal Z''(e)>0$ and $\mu$ has a continuous density positive at $e$. Then there exist $c_1,c_2>0$ and $T_1\ge2$ with $c_1T^{-1/2}\le\int(1+T\mathcal Z)^{-1}d\mu\le c_2T^{-1/2}$ for $T\ge T_1$, and for every continuous $\mathcal Q$ that is $C^1$ near $E$ and vanishes on $E$ there exists $K_2$ with
\[
\big|\int\mathcal Q(1+T\mathcal Z)^{-1}d\mu\big|\le K_2T^{-1}\log T\qquad(T\ge2).
\]
\end{lem}

\begin{proof}
For $e\in E$, Taylor's theorem yields positive constants $\delta_e,\alpha_e,\beta_e,\underline\varrho_e,\overline\varrho_e,\kappa_e$ with $\alpha_e(x-e)^2\le\mathcal Z(x)\le\beta_e(x-e)^2$, $\underline\varrho_e\le d\mu/dx\le\overline\varrho_e$, and $|\mathcal Q(x)|\le \kappa_e|x-e|$ for $|x-e|\le\delta_e$. Off $\bigcup_e(e-\delta_e,e+\delta_e)$, $\mathcal Z\ge z_0>0$, so those integrals are $O(T^{-1})$. For $T\ge\delta_e^{-2}$ (restrict to $|u|\le T^{-1/2}$ in the first bound),
\[
\frac{2}{(1+\beta_e)\sqrt T}\le\int_{-\delta_e}^{\delta_e}\frac{du}{1+\beta_eTu^2},\qquad \int_{-\delta_e}^{\delta_e}\frac{du}{1+\alpha_eTu^2}\le\frac{\pi}{\sqrt{\alpha_eT}},
\]
\[
\int_{-\delta_e}^{\delta_e}\frac{|u|\,du}{1+\alpha_eTu^2}=\frac{\log(1+\alpha_eT\delta_e^2)}{\alpha_eT}.
\]
The lemma follows.
\end{proof}

\begin{prop}\label{prop:fullrank}
Let $\nu_0,\dots,\nu_{S-1}>0$ and let $\Sigma_k$ be as in Lemma~\ref{lem:gapreg}. There exists $T_0>0$ such that for every $T\ge T_0$ the map $\zeta\to\Delta(\zeta,0)$ is differentiable at $0$ and the Jacobian matrix $\mathbf J=(\partial_{\zeta_b}\Delta_k(0,0))_{k,b}$ satisfies $\sum_{b\ne k}|\mathbf J_{kb}|<|\mathbf J_{kk}|$ for every $k$. In particular $\det\mathbf J\ne0$.
\end{prop}

\begin{proof}
As $\Delta$ is $C^\infty$ on $\mathcal P$, $\zeta\to\Delta(\zeta,0)$ is differentiable at $0$. Fix $k$, a regular $C^\infty$ parametrization $\alpha_k:[0,1]\to\overline{\Gamma_k}$ with $\alpha_k(0)=p_k$, and $\theta_k:=(F_k\circ\Sigma_k)'(0)\ne0$. If $F_b\ne F_k$, then $\mathbf J_{kb}=0$ (Lemma~\ref{lem:gapreg}). Let $F:=F_b=F_k$. Then $\nabla F\cdot X_0(0)=\ell_FF$ for a polynomial $\ell_F$, and $\nabla F\cdot Y_b=K_F^2P_b|\nabla F|^2$. Let $\gamma$ solve $X_0(0)$ on $\Gamma_k$ with $\gamma(0)=\Sigma_k(0)$, and $\mathcal E(t):=\exp(-\int_0^t\ell_F(\gamma))$. Dividing by $1+TZ\ge1$, let $\eta^u(\zeta,\cdot)$, $\eta^s(\zeta,\cdot)$ be the solutions of $X_0(0)+\sum_{b'}\zeta_{b'}Y_{b'}/(1+TZ)$ on the continued separatrices, with $\eta^u(\zeta,0)=\Sigma_k(u_k(\zeta,0))$ ($t\le0$) and $\eta^s(\zeta,0)=\Sigma_k(v_k(\zeta,0))$ ($t\ge0$). Then the functions
\[
z^{u,s}(t):=\partial_{\zeta_b}F(\eta^{u,s}(\zeta,t))|_{\zeta=0}
\]
satisfy
\[
\big(\mathcal E z^{u,s}\big)'=\mathcal E\,\frac{K_F^2P_b|\nabla F|^2}{1+TZ}\circ\gamma,\qquad z^u(0)-z^s(0)=\theta_k\mathbf J_{kb}.
\]
As $Y_{b'}$ has zero first jet at the corners, $z^u(t)=O(|\gamma(t)-p_k|^2)$ as $t\to-\infty$ and $z^s(t)=O(|\gamma(t)-p_{k+1}|^2)$ as $t\to\infty$. Differentiating the cofactor identity at $p_k,p_{k+1}$ and Proposition~\ref{prop:saddles}(2) with $T=0$ yield $\ell_F(p_k)=-\mu_k$ and $\ell_F(p_{k+1})=\lambda_{k+1}$, so $\mathcal E\asymp|\gamma-p_k|^{\nu_k}$, resp.\ $|\gamma-p_{k+1}|^{1/\nu_{k+1}}$, as $t\to-\infty$, resp.\ $\infty$, and the boundary terms vanish \cite[\S4.5]{GH83}, \cite[Ch.~8]{HY12}, and
\[
\mathbf J_{kb}=\theta_k^{-1}\int_0^1\frac{P_b(\alpha_k)\,d\mathfrak m_k}{1+TZ(\alpha_k)},
\]
where $\mathfrak m_k$, the image of $\mathcal E K_F(\gamma)^2|\nabla F(\gamma)|^2dt$ under $t\to\alpha_k^{-1}(\gamma(t))$, is independent of $T$ and has a continuous density, positive on $(0,1)$, $\asymp r^{\nu_k+1}$ near $0$ and $\asymp(1-r)^{1/\nu_{k+1}+1}$ near $1$ ($dt\asymp dr/r$, $K_F^2\asymp r^2$ near $0$). So $\mathfrak m_k$ is finite.

Let
\[
\mathcal M_k(T):=\int_0^1(1+TZ(\alpha_k))^{-1}d\mathfrak m_k.
\]
Lemma~\ref{lem:concentration} with $\mu=\mathfrak m_k$, $\mathcal Z=Z\circ\alpha_k$, $\mathcal Q=P_b\circ\alpha_k-\delta_{kb}$ ($F_b=F_k$) and Lemma~\ref{lem:samples} yield constants, uniform in $k,b$, with
\[
\Big|\frac{\theta_k\mathbf J_{kb}}{\mathcal M_k(T)}-\delta_{kb}\Big|\le\epsilon(T):=\frac{K_2\log T}{c_1\sqrt T}\qquad(T\ge T_1).
\]
Choose $T_0\ge\max(T_1,e^2)$ with $S\epsilon(T_0)<1$. The function $\epsilon$ decreases on $[e^2,\infty)$. For $T\ge T_0$ each row of $(\theta_k\mathbf J_{kb}/\mathcal M_k(T))$, hence of $\mathbf J$, is strictly diagonally dominant, and the Levy--Desplanques theorem implies $\det\mathbf J\ne0$. Order of choices: ratios, sections, $T_0$.
\end{proof}
\section{The host II: corner passages, calibration, limit cycles}\label{sec:hostB}

\subsection{Corner passages}\label{sec:dulac}

\begin{lem}\label{lem:dulac}
Let $\mathcal V\subset\Rr^d$ be open, $\mathcal K\subset\mathcal V$ compact, $r_0,c_0>0$, let $\mu,\lambda$ be positive $C^\infty$ functions on $\mathcal V$, and let $f_1,f_2$ be $C^\infty$ on $\mathcal V\times(-r_0,r_0)^2$ with $f_1(p,0,0)=f_2(p,0,0)=0$, $\mu-f_1\ge c_0$, and $\lambda+f_2\ge c_0$. Consider
\[
u'=u\big(-\mu(p)+f_1(p,u,v)\big),\qquad v'=v\big(\lambda(p)+f_2(p,u,v)\big).
\]
Then there exist $u_0,v_0\in(0,r_0)$, $t_0\in(0,v_0)$, $\eta\in(0,1)$ with $\eta<\min_{\mathcal K}\mu/\lambda$, $K_3>0$, a continuous positive function $\kappa$ on $\mathcal K$ and a function $\mathcal D$ on $\mathcal K\times(0,t_0)$ such that for $p\in\mathcal K$ and $t\in(0,t_0)$ the solution starting at $(u_0,t)$ reaches $\{v=v_0\}$ at a first time $\vartheta>0$, at the point $(\mathcal D(p,t),v_0)$, with $0<u\le u_0$ and $0<v<v_0$ on $[0,\vartheta)$, and
\[
\mathcal D(p,t)=\kappa(p)\,t^{\mu(p)/\lambda(p)}\big(1+r(p,t)\big),\qquad |r(p,t)|+t\,|\partial_tr(p,t)|\le K_3\,t^{\eta}.
\]
\end{lem}

\begin{proof}
Since $u'\le-c_0u$ and $v'\ge c_0v$ while $(u,v)\in(0,u_0]\times(0,v_0)$, the solution stays there until $v$ reaches $v_0$, which happens before the time $c_0^{-1}\log(v_0/t)$. Let $\nu:=\mu/\lambda$ and $\mathcal A:=(\mu-f_1)/(\lambda+f_2)$, so that $(u'/u)/(v'/v)=-\mathcal A$ and $\mathcal A(p,0,0)=\nu(p)$. The coordinates
\begin{gather*}
U:=u\,e^{G_1(p,u)},\qquad V:=v\,e^{G_2(p,v)},\\
G_1:=\int_0^u\Big(\frac{\nu}{\mathcal A(p,w,0)}-1\Big)\frac{dw}{w},\quad G_2:=\int_0^v\Big(\frac{\mathcal A(p,0,w)}{\nu}-1\Big)\frac{dw}{w},
\end{gather*}
are $C^\infty$, and a direct calculation shows $(U'/U)/(V'/V)=-\nu+E$ with
\[
E:=\nu-\nu^2\mathcal A(p,u,v)/(\mathcal A(p,u,0)\mathcal A(p,0,v)).
\]
Since $E=0$ on $\{u=0\}\cup\{v=0\}$, Hadamard's lemma implies $|E|+|U\partial_UE|\le M_EUV$ near $\mathcal K\times\{0\}$. Let $U_0,V_0$ be the values of $U,V$ at $u_0,v_0$, write $U=U_0\mathcal X$, $V=V_0y$, $\widetilde E(\mathcal X,y):=E(p,U_0\mathcal X,V_0y)$, and fix $0<\eta<c<\min_{\mathcal K}\nu$ with $\eta<1$. Taking $u_0,v_0$ small, $M:=M_E\max_{\mathcal K}(U_0V_0)$ satisfies $M\le\min_{\mathcal K}\nu-c$. The orbit through $(u_0,t)$ is the graph of the solution of
\[
y\,\partial_y\log\mathcal X(s,y)=-\nu+\widetilde E(\mathcal X,y),\qquad \mathcal X(s,s)=1,\qquad s:=V(t)/V_0,
\]
and integration (a first-exit argument) yields $(s/y)^{\nu+M}\le\mathcal X(s,y)\le(s/y)^c\le1$ on $[s,1]$. Thus $\log\mathcal X(s,1)=\nu\log s+I(s)$ with $I(s):=\int_s^1\widetilde E(\mathcal X,y)\,dy/y$ and $|I(s)|\le M\int_s^1(s/y)^\eta dy\le Ms^\eta/(1-\eta)$. The function $\mathcal Y:=s\,\partial_s\log\mathcal X$ satisfies $\partial_y\mathcal Y=\mathcal B\mathcal Y$ with
\[
\mathcal B:=\mathcal X\,\partial_{\mathcal X}\widetilde E(\mathcal X,y)/y,\qquad |\mathcal B|\le M\mathcal X\le M(s/y)^\eta,
\]
and $\mathcal Y(s,s)=\nu-\widetilde E(1,s)$ (differentiate $\mathcal X(s,s)=1$). Hence $|\mathcal Y|\le(\nu+M)e^{M/(1-\eta)}$ and
\[
s\,I'(s)=\mathcal Y(s,1)-\nu=-\widetilde E(1,s)+\int_s^1\mathcal B\mathcal Y\,dy=O(s^\eta).
\]
Returning to $u$ through $s=te^{G_2(p,t)}/V_0$ and $\mathcal D=U_0\mathcal X(s,1)e^{-G_1(p,\mathcal D)}$, the lemma holds with $\kappa:=U_0V_0^{-\nu}$ and $1+r=e^{I(s)+\nu G_2(p,t)-G_1(p,\mathcal D)}$, because $\eta<\min(1,\nu)$. All constants are uniform on $\mathcal K$ by compactness. This is the classical finitely smooth Dulac asymptotics \cite[\S24]{IY08}, \cite[Ch.~5]{Rou98}, \cite{Mou90}.
\end{proof}

\begin{prop}\label{prop:passages}
Let $\nu_0,\dots,\nu_{S-1}>0$, let $\Sigma_k,u_k,v_k,\Delta_k$ be as in Lemma~\ref{lem:gapreg}, and let $T>0$. Then there exist an open set $\mathcal P\ni0$ as in Lemma~\ref{lem:gapreg}, a closed ball $\mathcal P_0\subset\mathcal P$ centred at $0$, numbers $0<t_*<\epsilon<1$, $C_1>0$, $\eta\in(0,1)$ with $\eta<\nu_k(\tau)$ for $(\zeta,\tau)\in\mathcal P_0$ and $k\in\ZZ/S$, and for each $k$ functions $R_k,\vartheta_k\in C^\infty(\mathcal P\times(0,\epsilon))$ and a continuous positive function $\kappa_k$ on $\mathcal P$ such that:
\begin{enumerate}
\item For $p=(\zeta,\tau)\in\mathcal P$ and $t\in(0,\epsilon)$ the solution of $X_{T,\zeta}(\tau)$ starting at $\Sigma_{k-1}(v_{k-1}(p)+t)$ is at $\Sigma_k(v_k(p)+R_k(p,t))$ at the time $\vartheta_k(p,t)>0$.
\item For $p=(\zeta,\tau)\in\mathcal P_0$ and $t\in(0,t_*)$,
\[
R_k(p,t)=\Delta_k(p)+\kappa_k(p)\,t^{\nu_k(\tau)}\big(1+r_k(p,t)\big),\qquad |r_k(p,t)|+t\,|\partial_tr_k(p,t)|\le C_1t^{\eta}.
\]
\item If $p\in\mathcal P$ and $t_0,\dots,t_S$ satisfy $t_k\in(0,\epsilon)$ and $t_{k+1}=R_k(p,t_k)$ for $0\le k\le S-1$ (where $\Sigma_{-1}:=\Sigma_{S-1}$, $v_{-1}:=v_{S-1}$), then the $S$ solutions of (1), traversed one after another, form a solution of duration $\sum_{k=0}^{S-1}\vartheta_k(p,t_k)$ from $\Sigma_{S-1}(v_{S-1}(p)+t_0)$ to $\Sigma_{S-1}(v_{S-1}(p)+t_S)$ whose points at intermediate times are not in $\Sigma_{S-1}((-r_{S-1},r_{S-1}))$.
\end{enumerate}
\end{prop}

\begin{proof}
Fix $k$. The stable and unstable manifold theorem with parameters \cite[Ch.~IX]{Har02} and Hadamard's lemma yield $C^\infty$ coordinates $(u,v)$ near $p_k$, depending smoothly on $p$ near $0$, in which the stable and unstable separatrices of $p_k$ are $\{v=0\}$ and $\{u=0\}$, the quadrant $\{u,v>0\}$ faces $\Omega$, and $X_{T,\zeta}(\tau)$ has the form of Lemma~\ref{lem:dulac} with $\mu/\lambda=\mu_k/\lambda_k=\nu_k(\tau)$ (Proposition~\ref{prop:saddles}(2)). The inequalities $\mu-f_1\ge c_0$ and $\lambda+f_2\ge c_0$ hold after shrinking. We use Lemma~\ref{lem:dulac} on a closed ball $\mathcal K$ centred at $0$, one for all $k$, with output exponent $\eta_k$, and let $\mathcal P$ be the interior of $\mathcal K$, intersected with the set of Lemma~\ref{lem:gapreg}. By the flow box theorem with parameters \cite[Ch.~2]{PdM82}, \cite[Ch.~V]{Har02}, the transitions from $\Sigma_{k-1}$ (coordinate $t$ centred at $v_{k-1}(p)$) to $\{u=u_0\}$, and from $\{v=v_0\}$ to $\Sigma_k$ (coordinate centred at $v_k(p)$), are $C^\infty$ maps $h(p,\cdot)$, $j(p,\cdot)$ with $h'(p,0)>0$, $j'(p,0)>0$ ($\Sigma_{k-1}$ and $\Sigma_k$ point into $\Omega$), $h(p,0)=0$ (the stable separatrix of $p_k$ passes through $\Sigma_{k-1}(v_{k-1}(p))$ and $(u_0,0)$) and $j(p,0)=u_k(p)-v_k(p)=\Delta_k(p)$ (the unstable separatrix passes through $(0,v_0)$ and $\Sigma_k(u_k(p))$). Thus $R_k(p,t)=j(p,\mathcal D(p,h(p,t)))$, and Taylor expansion of $h,j$ with $\mathcal D=O(t^{\nu_k(\tau)})$ and $\eta_k<\min(1,\nu_k(\tau))$ implies (2) with $\eta_k$ in place of $\eta$ and $\kappa_k(p):=j'(p,0)\,\kappa(p)\,h'(p,0)^{\nu_k(\tau)}$. Smoothness of $R_k,\vartheta_k$ for $t>0$ follows from smooth dependence on initial data at transverse crossings \cite[Ch.~V]{Har02}. Let $\eta:=\min_k\eta_k$, and $\mathcal P_0$, $t_*$, $\epsilon$ small. Item (3) holds because for $\epsilon$ and $\mathcal P$ small the solutions of (1) stay in a thin tube around $\Gamma$ which they traverse once, crossing $\Sigma_{S-1}$ once per turn.
\end{proof}

\subsection{Calibration and the rescaled return map}\label{sec:calibration}

\begin{lem}\label{lem:calibration}
Let $\varkappa>0$ and let $K$ be a real function continuous at $0$ with $K(0)>0$. Let $\tau(\sigma):=-\log K(0)/(\varkappa\log\sigma)$ for $\sigma\in(0,1)$. Then, as $\sigma\to0^+$, $\tau(\sigma)\to0$ and
\[
K(\tau(\sigma))\,\sigma^{\varkappa\tau(\sigma)}=K(\tau(\sigma))/K(0)\to1.
\]
\end{lem}

\begin{proof}
As $\sigma\to0^+$, $\log\sigma\to-\infty$, so $\tau(\sigma)\to0$ and $K(\tau(\sigma))\to K(0)$. Moreover $\sigma^{\varkappa\tau(\sigma)}=e^{-\log K(0)}=1/K(0)$.
\end{proof}

\begin{lem}\label{lem:gapinverse}
Let $E$ be a finite-dimensional normed space and $\Delta:E\times\Rr\to E$ a map that is $C^1$ on an open set containing $(0,0)$, with $\Delta(0,\tau)=0$ for $|\tau|$ small and $D_\zeta\Delta(0,0)$ invertible. Then there exist an open neighbourhood $U_3$ of $0$ in $E$, $\epsilon_3>0$, $C_3>0$ and a $C^1$ map $\mathcal Q$ on $U_3\times(-\epsilon_3,\epsilon_3)$ with $\Delta(\mathcal Q(d,\tau),\tau)=d$ and $\|\mathcal Q(d,\tau)\|\le C_3\|d\|$.
\end{lem}

\begin{proof}
The implicit function theorem \cite[Thm.~9.28]{Rud76} for $\Delta(\zeta,\tau)-d=0$ yields $\mathcal Q$ with $\mathcal Q(0,\tau)=0$. On a convex $U_3$ with $\|D_d\mathcal Q\|\le C_3$ the mean value inequality implies the bound.
\end{proof}

\begin{prop}\label{prop:return}
Let $s\ge2$, $0<x_-\le x_+$, and let $W$ be a word with $|W|\le S$ that is legal on $[x_-,x_+]$ with $\Pi(W)=1$. Let $(g_0,\dots,g_{S-1})$, $g_k(x)=a_k+c_kx^{\nu_k}$, be $W$ followed by $S-|W|$ identity segments, and let $X_{T,\zeta}(\tau)$, $\Sigma_k$, $v_k$ be the objects of Sections~\ref{sec:field}--\ref{sec:fullrank} for these $\nu_k$. Then there exist $T>0$, $\epsilon\in(0,1)$, maps $\sigma\to\zeta(\sigma)$, $\sigma\to\tau(\sigma)$ on $(0,\epsilon)$ tending to $0$ as $\sigma\to0^+$, an open set $U\supset[x_-,x_+]$, and for $\sigma\in(0,\epsilon)$ functions $\mathbf P_\sigma,\vartheta_\sigma\in C^\infty(U)$ with $\vartheta_\sigma>0$, such that $(\zeta(\sigma),\tau(\sigma))$ is in the set $\mathcal P$ of Lemma~\ref{lem:gapreg} for $T$ and, with $X_\sigma:=X_{T,\zeta(\sigma)}(\tau(\sigma))$ and $o_\sigma:=v_{S-1}(\zeta(\sigma),\tau(\sigma))$:
\begin{enumerate}
\item $\deg X_\sigma\le20s$.
\item For $x\in U$ the solution of $X_\sigma$ starting at $\Sigma_{S-1}(o_\sigma+\sigma x)$ is at $\Sigma_{S-1}(o_\sigma+\sigma\mathbf P_\sigma(x))$ at the time $\vartheta_\sigma(x)$ and does not meet $\Sigma_{S-1}((-r_{S-1},r_{S-1}))$ at times in $(0,\vartheta_\sigma(x))$.
\item $\mathbf P_\sigma\to W$ and $\mathbf P_\sigma'\to W'$ uniformly on $[x_-,x_+]$ as $\sigma\to0^+$.
\end{enumerate}
\end{prop}

\begin{proof}
The identity segment is $x\to x$, so $G_k:=g_{k-1}\circ\dots\circ g_0$ satisfies $G_S=W$ and $G_k>0$ on $[x_-,x_+]$. Let $\pi_k:=\nu_0\cdots\nu_{k-1}$ $(0\le k\le S)$. Then $\pi_0=\pi_S=1$, $\pi_{k+1}=\nu_k\pi_k$, and $\pi_{\min}:=\min_k\pi_k>0$. By compactness there exist $0<l_k^-<l_k^+$ with $G_k([x_-,x_+])\subset(l_k^-,l_k^+)$.

\emph{Choice of $T$ and of the scales.} Let $T:=T_0$ (Proposition~\ref{prop:fullrank}), take the data of Proposition~\ref{prop:passages} for this $T$ (write $\epsilon_1$ for the number $\epsilon$ there), and let $\delta_1>0$ with $(0,\tau)\in\mathcal P_0$ for $|\tau|<\delta_1$. On $(-\delta_1,\delta_1)$ let $\rho_0:=1$, $\rho_{k+1}(\tau):=\kappa_k(0,\tau)\rho_k(\tau)^{\nu_k}/c_k$ for $k+1<S$, $\rho_S:=1$. The functions $\rho_k$ are continuous and positive. Lemma~\ref{lem:calibration} with $\varkappa:=\pi_{S-1}$ and
\[
K(\tau):=\kappa_{S-1}(0,\tau)\,\rho_{S-1}(\tau)^{\nu_{S-1}+\tau}/c_{S-1}
\]
(continuous, with $K(0)>0$) defines $\tau(\sigma)$. Let $\varsigma_k(\sigma):=\rho_k(\tau(\sigma))\sigma^{\pi_k}$, so $\varsigma_0=\varsigma_S=\sigma$ and $\varsigma_k\to0$.

\emph{Choice of $\zeta$.} By Lemma~\ref{lem:gapreg} and Proposition~\ref{prop:fullrank}, the map $\Delta=(\Delta_k)_k$ on $\mathcal P$ satisfies the hypotheses of Lemma~\ref{lem:gapinverse} with $E=\Rr^S$. Let $\zeta(\sigma):=\mathcal Q(d(\sigma),\tau(\sigma))$ with
\[
d_k(\sigma):=\varsigma_{k+1}(\sigma)a_k=O(\sigma^{\pi_{\min}}).
\]
Thus, for $\sigma$ small, $p(\sigma):=(\zeta(\sigma),\tau(\sigma))\in\mathcal P_0$, $\zeta(\sigma)=O(\sigma^{\pi_{\min}})$ and
\[
\Delta_k(p(\sigma))=\varsigma_{k+1}(\sigma)\,a_k\qquad(k\in\ZZ/S).
\]

\emph{Normalized passages.} For $x$ near $[l_k^-,l_k^+]$ and $\sigma$ small, $\varsigma_kx<t_*$, and Proposition~\ref{prop:passages}(2) implies
\[
\widehat R_k(x):=\frac{R_k(p(\sigma),\varsigma_kx)}{\varsigma_{k+1}}=a_k+m_k(\sigma)\,x^{\nu_k(\tau(\sigma))}\big(1+r_k(p(\sigma),\varsigma_kx)\big),\qquad m_k:=\frac{\kappa_k(p(\sigma))\,\varsigma_k^{\nu_k(\tau(\sigma))}}{\varsigma_{k+1}}.
\]
For $k<S-1$ the recursion implies $m_k=c_k\kappa_k(p(\sigma))/\kappa_k(0,\tau(\sigma))$. For $k=S-1$, $\pi_{S-1}\nu_{S-1}=1$ implies
\[
m_{S-1}=c_{S-1}\,\kappa_{S-1}(p(\sigma))\kappa_{S-1}(0,\tau(\sigma))^{-1}K(\tau(\sigma))\sigma^{\pi_{S-1}\tau(\sigma)}.
\]
By continuity of $\kappa_k$ and Lemma~\ref{lem:calibration}, $m_k\to c_k$. The remainder satisfies $|r_k|+x|\partial_xr_k|\le C_1(\varsigma_kx)^\eta\to0$, and $x^{\nu_k(\tau(\sigma))}\to x^{\nu_k}$ in $C^1$ on $[l_k^-,l_k^+]$. Hence $\widehat R_k\to g_k$ and $\widehat R_k'\to g_k'$ uniformly on $[l_k^-,l_k^+]$.

\emph{Composition.} Since $G_k$ is $C^1$ near $[x_-,x_+]$ (Lemma~\ref{lem:legal}(1)) and $G_k([x_-,x_+])\subset(l_k^-,l_k^+)$, induction on $k$ (uniform continuity and the chain rule) yields $\epsilon\in(0,1)$ and a bounded open $U\supset[x_-,x_+]$ such that for $\sigma\in(0,\epsilon)$ the maps $\widehat R_{k-1}\circ\dots\circ\widehat R_0$ are $C^\infty$ on $U$ with values in $(l_k^-,l_k^+)$ and converge to $G_k$ in $C^1$ on $\overline U$. Let $\mathbf P_\sigma:=\widehat R_{S-1}\circ\dots\circ\widehat R_0$. Then (3) holds. For $x\in U$ let
\[
t_k:=\varsigma_k\,\widehat R_{k-1}\circ\dots\circ\widehat R_0(x)\in(0,\epsilon_1),
\]
where the inclusion holds as $\varsigma_k\to0$. Then $t_0=\sigma x$, $t_{k+1}=R_k(p(\sigma),t_k)$, and $t_S=\sigma\mathbf P_\sigma(x)$ because $\varsigma_0=\varsigma_S=\sigma$. Proposition~\ref{prop:passages}(3) implies (2) with $\vartheta_\sigma(x):=\sum_k\vartheta_k(p(\sigma),t_k)$, and Proposition~\ref{prop:saddles}(1) implies (1).
\end{proof}

\subsection{From fixed points to hyperbolic limit cycles}\label{sec:cycles}

\begin{lem}\label{lem:persist}
Let $0<x_-\le x_+$, $\sigma_0>0$, let $W$ be $C^1$ near $[x_-,x_+]$, and for $\sigma\in(0,\sigma_0)$ let $\mathbf P_\sigma$ be $C^1$ near $[x_-,x_+]$ with $\mathbf P_\sigma\to W$ and $\mathbf P_\sigma'\to W'$ uniformly on $[x_-,x_+]$ as $\sigma\to0^+$. Let $p_1,\dots,p_M\in(x_-,x_+)$ be pairwise distinct with $W(p_r)=p_r$ and $W'(p_r)\ne1$. Then there exist $\sigma_1\in(0,\sigma_0]$ and pairwise disjoint intervals $[l_r,u_r]\subset(x_-,x_+)$ such that for every $\sigma\in(0,\sigma_1)$ and every $r$ the map $\mathbf P_\sigma$ has exactly one fixed point $x_r$ in $[l_r,u_r]$, and $\mathbf P_\sigma'(x_r)\ne1$.
\end{lem}

\begin{proof}
If $M=0$ let $\sigma_1:=\sigma_0$. Otherwise $x_-<x_+$. Let $\epsilon$ and $[l_r,u_r]$ be given by Lemma~\ref{lem:stability} for $f:=W-\id$, and $\sigma_1$ such that for $\sigma<\sigma_1$
\[
\sup_{[x_-,x_+]}(|\mathbf P_\sigma-W|+|\mathbf P_\sigma'-W'|)<\epsilon.
\]
Lemma~\ref{lem:stability} for $g:=\mathbf P_\sigma-\id$ implies the claim.
\end{proof}

\begin{lem}\label{lem:multiplier}
Let $X$ be a polynomial vector field, $I_0\subset\Rr$ an open interval, $\iota:I_0\to\Rr^2$ a $C^\infty$ embedding with $\det(X(\iota(x)),\iota'(x))\ne0$, $V_0\subset I_0$ an open interval, $\mathbf P,\vartheta\in C^\infty(V_0)$ with $\mathbf P(V_0)\subset I_0$ and $\vartheta>0$, and $\Xi$ a $C^\infty$ map on an open set containing $\{(t,x):x\in V_0,\ 0\le t\le\vartheta(x)\}$ with $\partial_t\Xi=X(\Xi)$, $\Xi(0,x)=\iota(x)$, and $\Xi(\vartheta(x),x)=\iota(\mathbf P(x))$. If $a\in V_0$ and $\mathbf P(a)=a$, then
\[
\mathbf P'(a)=\exp\int_0^{\vartheta(a)}\operatorname{div}X(\Xi(t,a))\,dt .
\]
\end{lem}

\begin{proof}
Let $\vartheta:=\vartheta(a)$, $w(t):=\partial_x\Xi(t,a)$, and $d(t):=\det(X(\Xi(t,a)),w(t))$. The vectors $X(\Xi(t,a))$ and $w(t)$ solve the variational equation $\dot y=DX(\Xi(t,a))y$, so $d'=(\operatorname{div}X)\,d$ (Liouville) and
\[
d(\vartheta)=d(0)\exp\int_0^{\vartheta}\operatorname{div}X(\Xi(t,a))\,dt.
\]
Differentiating $\Xi(\vartheta(x),x)=\iota(\mathbf P(x))$ at $a$ yields
\[
\vartheta'(a)X(\iota(a))+w(\vartheta)=\mathbf P'(a)\iota'(a).
\]
Applying $\det(X(\iota(a)),\cdot)$ and using $\Xi(\vartheta,a)=\iota(a)=\Xi(0,a)$ yields $d(\vartheta)=\mathbf P'(a)\,d(0)$, and $d(0)=\det(X(\iota(a)),\iota'(a))\ne0$. The lemma follows (cf.\ \cite[\S3.4]{Per01}).
\end{proof}

\begin{lem}\label{lem:isolation}
In the setting of Lemma~\ref{lem:multiplier} assume moreover $\Xi(t,x)\notin\iota(I_0)$ for $0<t<\vartheta(x)$. If $a\in V_0$, $\mathbf P(a)=a$ and $\mathbf P'(a)\ne1$, then $O(a):=\Xi([0,\vartheta(a)],a)$ is a hyperbolic limit cycle of $X$ and $O(a)\cap\iota(I_0)=\{\iota(a)\}$.
\end{lem}

\begin{proof}
The $\vartheta(a)$-periodic extension of $\Xi(\cdot,a)$ is a periodic solution, not constant since $X(\iota(a))\ne0$, with image $O(a)$. The hypothesis on $\Xi$ implies $O(a)\cap\iota(I_0)=\{\iota(a)\}$, and Lemma~\ref{lem:multiplier} implies
\[
\int_0^{\vartheta(a)}\operatorname{div}X(\Xi(t,a))\,dt=\log\mathbf P'(a)\ne0.
\]
Isolation: since $\mathbf P'(a)>0$ and $\mathbf P'(a)\ne1$, there exists an open interval $V_1\subset V_0$ with $a\in V_1$ on which $\mathbf P$ is strictly increasing and $a$ is the only fixed point, so any $x_0,\dots,x_m\in V_1$ with $x_{i+1}=\mathbf P(x_i)$ and $x_m=x_0$ equal $a$. By the flow box theorem along $O(a)$ \cite[Ch.~V]{Har02} there exists an open $\mathcal U\supset O(a)$ such that any periodic orbit $O'\subset\mathcal U$ crosses $\iota(V_1)$, and by uniqueness of solutions the points of $O'\cap\iota(V_1)$, in the order of time over one period, are $\iota(x_0),\dots,\iota(x_m)$ for such a closed chain. Thus $\iota(a)\in O'$, and $O'=O(a)$ by uniqueness of solutions.
\end{proof}

\begin{proof}[Proof of Theorem~\ref{thm:host}]
If $M=0$ take $X:=0$. Let $M\ge1$. Since $\sqrt{|W|}/2\le s$, $|W|\le4s^2<S$. Let $T,\epsilon,\zeta(\sigma),\tau(\sigma),U,\mathbf P_\sigma,\vartheta_\sigma$ be given by Proposition~\ref{prop:return}. By Lemma~\ref{lem:legal}(1), $W$ is $C^1$ near $[x_-,x_+]$, so Lemma~\ref{lem:persist} with $\sigma_0:=\epsilon$ yields $\sigma_1$ and the intervals $[l_r,u_r]$. Fix $\sigma:=\sigma_1/2$ and $X:=X_\sigma$. Then $\deg X\le20s$, and $\mathbf P_\sigma$ has fixed points $x_r\in[l_r,u_r]$, pairwise distinct, with $\mathbf P_\sigma'(x_r)\ne1$. Let $I_0:=\{x:|o_\sigma+\sigma x|<r_{S-1}\}$ and $\iota(x):=\Sigma_{S-1}(o_\sigma+\sigma x)$, transverse to $X$ by Lemma~\ref{lem:gapreg}. Fix $r$, and let $V_0\subset U$ be an open interval with $x_r\in V_0$, and $\Xi(t,x)$ the solution of $X$ at time $t$ starting at $\iota(x)$, which is $C^\infty$ on an open set \cite[Ch.~V]{Har02}. Proposition~\ref{prop:return}(2) implies $V_0\cup\mathbf P_\sigma(V_0)\subset I_0$ and the remaining hypotheses of Lemma~\ref{lem:isolation}. Hence $O_r:=O(x_r)$ is a hyperbolic limit cycle of $X$ with $O_r\cap\iota(I_0)=\{\iota(x_r)\}$. As $\iota$ is injective, the $O_r$ are pairwise distinct. The theorem follows.
\end{proof}
\appendix

\section{The data}\label{app:data}

The numbers $a,\gamma,\alpha_i,\beta_i$ $(1\le i\le49)$ of Theorem~\ref{thm:cert} and $u_k,v_k$ $(1\le k\le183)$ of Proposition~\ref{prop:checks} are dyadic rationals, with
\begin{equation}\label{eq:a-gamma}
a=3619630669166631734475479811931096\cdot2^{-108}\approx11.1538481,\qquad\gamma=2^{-360}.
\end{equation}
Table~\ref{tab:pairs} lists integers $m_i,e_i,n_i,f_i$ and Table~\ref{tab:brackets} lists integers $t_k$ such that
\[
\alpha_i=m_i\,2^{e_i},\quad\beta_i=n_i\,2^{f_i}\quad(1\le i\le49),\qquad
u_k=t_k\,2^{-82},\quad v_k=u_k+2^{-82}\quad(1\le k\le183).
\]
Tables~\ref{tab:pairs} and~\ref{tab:brackets} are generated mechanically from the data file \nolinkurl{hilbert16_367_roots.json}. The data file and the script \nolinkurl{verify_template.py} of Section~\ref{sec:verification} accompany the paper as ancillary files.

\begin{center}
\small
\begin{longtable}{@{}rrrrr@{}}
\caption{The coefficients $\alpha_i=m_i\cdot 2^{e_i}$ and $\beta_i=n_i\cdot 2^{f_i}$ of the template, $1\le i\le 49$.}\label{tab:pairs}\\
\hline
$i$ & $m_i$ & $e_i$ & $n_i$ & $f_i$\\
\hline
\endfirsthead
\hline
$i$ & $m_i$ & $e_i$ & $n_i$ & $f_i$\\
\hline
\endhead
\hline
\endfoot
1 & \texttt{3105261365289511707859521755268099} & $179$ & \texttt{2706267646357366187558280189123413} & $-83$\\
2 & \texttt{4002975836278912222185099759295522} & $167$ & \texttt{2626767494772052301221389783699045} & $-84$\\
3 & \texttt{5160201066933202381125296033257751} & $155$ & \texttt{5099204478959908258597472465460484} & $-86$\\
4 & \texttt{3325985018336835873644813093422927} & $144$ & \texttt{4949407659182537413122112411669380} & $-87$\\
5 & \texttt{4287498182899133193121482013361734} & $132$ & \texttt{4804011346939007284296035516267903} & $-88$\\
6 & \texttt{2763488194938453081136149784903066} & $121$ & \texttt{4662886270597748793260487840703470} & $-89$\\
7 & \texttt{3562388450275082896125813345802967} & $109$ & \texttt{4525906956153098142687444811608463} & $-90$\\
8 & \texttt{4592243778703389538213891297585819} & $97$ & \texttt{4392951615596255917494968793972514} & $-91$\\
9 & \texttt{2959910635442416713725110863915239} & $86$ & \texttt{4263902038626052793724310921340363} & $-92$\\
10 & \texttt{3815594899572999328842059835085608} & $74$ & \texttt{4138643487547710743851204590912351} & $-93$\\
11 & \texttt{4918649996833905719010134531986307} & $62$ & \texttt{4017064595259865035029202185774592} & $-94$\\
12 & \texttt{3170294335237424753538660730996468} & $51$ & \texttt{3899057266237715858816256624007040} & $-95$\\
13 & \texttt{4086798686024864505915751404759571} & $39$ & \texttt{3784516580425480241857439196307624} & $-96$\\
14 & \texttt{2634128212395734002561484903029333} & $28$ & \texttt{3673340699952202015314274116156778} & $-97$\\
15 & \texttt{3395631628769164086975904918335742} & $16$ & \texttt{3565430778587925511203502023904635} & $-98$\\
16 & \texttt{4377279019311945336571619058988786} & $4$ & \texttt{3460690873859730669515390993959557} & $-99$\\
17 & \texttt{2821356040297839755728902088281505} & $-7$ & \texttt{3359027861749492116891533222950669} & $-100$\\
18 & \texttt{3636985383388389970431929824238147} & $-19$ & \texttt{3260351353897520182727129568154228} & $-101$\\
19 & \texttt{4688406032435523547897131539997923} & $-31$ & \texttt{3164573617238469792662667466732373} & $-102$\\
20 & \texttt{3021891595354600124266026620084506} & $-42$ & \texttt{3071609495998065703431495423244425} & $-103$\\
21 & \texttt{3895494012633964040911154315954175} & $-54$ & \texttt{2981376335981291539628406180715603} & $-104$\\
22 & \texttt{5021647244326905002007158061164790} & $-66$ & \texttt{2893793911084727428373204076791156} & $-105$\\
23 & \texttt{3236680760472456389539942843219116} & $-77$ & \texttt{2808784351967698513452765406567882} & $-106$\\
24 & \texttt{4172376845883670430738085908056890} & $-89$ & \texttt{2726272076818816024513394851410291} & $-107$\\
25 & \texttt{2689287240908649920864505462661660} & $-100$ & \texttt{2646183724156355585771294940952526} & $-108$\\
26 & \texttt{3466736649758362151090766550756811} & $-112$ & \texttt{5136896175205451457819593491158103} & $-110$\\
27 & \texttt{4468939879667569857048299046981096} & $-124$ & \texttt{4985992105150069468151168208428658} & $-111$\\
28 & \texttt{2880435646802770308647239049287632} & $-135$ & \texttt{4839521069670935043361380771636872} & $-112$\\
29 & \texttt{3713144387160238235969417863502448} & $-147$ & \texttt{4697352841693675741456968659204389} & $-113$\\
30 & \texttt{4786581937771594749393934033389301} & $-159$ & \texttt{4559361019760572433561045735887719} & $-114$\\
31 & \texttt{3085170445596861516568679993367340} & $-170$ & \texttt{4425422915647302222897793714313603} & $-115$\\
32 & \texttt{3977066224762296527338485450033548} & $-182$ & \texttt{4295419445281108980995298979649254} & $-116$\\
33 & \texttt{5126801269187263765262900585804190} & $-194$ & \texttt{4169235022863417087281675984319172} & $-117$\\
34 & \texttt{3304457327123223921629333039444672} & $-205$ & \texttt{4046757458102753021629319938493175} & $-118$\\
35 & \texttt{4259747024877126430672133058056595} & $-217$ & \texttt{3927877856466604827249515706086535} & $-119$\\
36 & \texttt{2745601307514309816587169644606069} & $-228$ & \texttt{3812490522363533592337146409788197} & $-120$\\
37 & \texttt{3539330620245803329921615994772987} & $-240$ & \texttt{3700492865169456379467840074988988} & $-121$\\
38 & \texttt{4562520131792388461330988969158408} & $-252$ & \texttt{3591785308014548778306799427760990} & $-122$\\
39 & \texttt{2940752388874755689692626119203138} & $-263$ & \texttt{3486271199249669718002761327462273} & $-123$\\
40 & \texttt{3790898171566048395561534571173359} & $-275$ & \texttt{3383856726513589434439838164269515} & $-124$\\
41 & \texttt{4886813660865986311004444108365144} & $-287$ & \texttt{3284450833324447540431274233919531} & $-125$\\
42 & \texttt{3149774364171486965988226273006009} & $-298$ & \texttt{3187965138121307880727655117044478} & $-126$\\
43 & \texttt{4060346570847419746884240033678767} & $-310$ & \texttt{3094313855684568133826732842900892} & $-127$\\
44 & \texttt{2617078617231767441843122750909479} & $-321$ & \texttt{3003413720862772136462731352865132} & $-128$\\
45 & \texttt{3373653144509271977521898649164760} & $-333$ & \texttt{2915183914521895197101006208280783} & $-129$\\
46 & \texttt{4348946741745907195207519894292168} & $-345$ & \texttt{2829545991153044378062622963274820} & $-130$\\
47 & \texttt{2803094426282263949564668835020870} & $-356$ & \texttt{2746423777614119357753637705243448} & $-131$\\
48 & \texttt{3613426389728961460623640147506083} & $-368$ & \texttt{2665739081601206220427537379440381} & $-132$\\
49 & \texttt{4645762021773806302751483938773761} & $-380$ & \texttt{5168941948639770331451724648771617} & $-134$\\
\end{longtable}
\end{center}

\begin{center}
\small
\begin{longtable}{@{}rl@{\qquad}rl@{}}
\caption{The left endpoints $u_k=t_k\cdot 2^{-82}$ of the brackets $[u_k,u_k+2^{-82}]$, $1\le k\le 183$.}\label{tab:brackets}\\
\hline
$k$ & $t_k$ & $k$ & $t_k$\\
\hline
\endfirsthead
\hline
$k$ & $t_k$ & $k$ & $t_k$\\
\hline
\endhead
\hline
\endfoot
1 & \texttt{437006100361935730564502} & 93 & \texttt{40641446387545731197509456}\\
2 & \texttt{874011447011279962589456} & 94 & \texttt{41078450668229286167406528}\\
3 & \texttt{1311015727758062868929567} & 95 & \texttt{41515456015037609535622410}\\
4 & \texttt{1748019254803578100650230} & 96 & \texttt{41952462115003123280745728}\\
5 & \texttt{2185022781849093332393289} & 97 & \texttt{42389468216344547698166616}\\
6 & \texttt{2622027062595876238713423} & 98 & \texttt{42826473560598004624908372}\\
7 & \texttt{3059032409245220470753974} & 99 & \texttt{43263477847140656991453274}\\
8 & \texttt{3496038509607156201308400} & 100 & \texttt{43700481360335594283923361}\\
9 & \texttt{3933044609969091931877304} & 101 & \texttt{44137484920037602282300724}\\
10 & \texttt{4370049956618436163902785} & 102 & \texttt{44574489124924414282579431}\\
11 & \texttt{4807054237365219070238794} & 103 & \texttt{45011494644918570539731441}\\
12 & \texttt{5244057764410734301965550} & 104 & \texttt{45448500356333419228634609}\\
13 & \texttt{5681061291456249533702029} & 105 & \texttt{45885507311918668526375349}\\
14 & \texttt{6118065572203032440027972} & 106 & \texttt{46322510819974298501716581}\\
15 & \texttt{6555070918852376672064431} & 107 & \texttt{46759518955323772135426018}\\
16 & \texttt{6992077019214312402620642} & 108 & \texttt{47196514625583255362713270}\\
17 & \texttt{7429083119576248133190255} & 109 & \texttt{47633533667894347999068351}\\
18 & \texttt{7866088466225592365212770} & 110 & \texttt{48070508390575746380890968}\\
19 & \texttt{8303092746972375271553426} & 111 & \texttt{48507567790077554560706356}\\
20 & \texttt{8740096274017890503274552} & 112 & \texttt{48944479618063666774177680}\\
21 & \texttt{9177099801063405735017100} & 113 & \texttt{49381641291417472059429172}\\
22 & \texttt{9614104081810188641336716} & 114 & \texttt{49818406393333046977273421}\\
23 & \texttt{10051109428459532873379988} & 115 & \texttt{50255752615779053716707796}\\
24 & \texttt{10488115528821468603928397} & 116 & \texttt{50692318242880119799367636}\\
25 & \texttt{10925121629183404334507321} & 117 & \texttt{51129805437542547619180295}\\
26 & \texttt{11362126975832748566518742} & 118 & \texttt{51566390216613649344745155}\\
27 & \texttt{11799131256579531472872123} & 119 & \texttt{52003595741879344387114327}\\
28 & \texttt{12236134783625046704579521} & 120 & \texttt{52440749016261850879005317}\\
29 & \texttt{12673138310670561936335807} & 121 & \texttt{52877285660832874387112355}\\
30 & \texttt{13110142591417344842642783} & 122 & \texttt{53314811889378763032230430}\\
31 & \texttt{13547147938066689074696686} & 123 & \texttt{53751658977922650160918530}\\
32 & \texttt{13984154038428624805236939} & 124 & \texttt{54188257803734241884727802}\\
33 & \texttt{14421160138790560535821599} & 125 & \texttt{54625872935020606110195844}\\
34 & \texttt{14858165485439904767829261} & 126 & \texttt{55062748806630186941004597}\\
35 & \texttt{15295169766186687674185018} & 127 & \texttt{55499187440584898919340556}\\
36 & \texttt{15732173293232202905890861} & 128 & \texttt{55936737309983189413089630}\\
37 & \texttt{16169176820277718137648365} & 129 & \texttt{56374034833189489208566784}\\
38 & \texttt{16606181101024501043954002} & 130 & \texttt{56810276687683765667425374}\\
39 & \texttt{17043186447673845276009835} & 131 & \texttt{57247314275986766498301492}\\
40 & \texttt{17480192548035781006547128} & 132 & \texttt{57685164202856605755005539}\\
41 & \texttt{17917198648397716737136054} & 133 & \texttt{58121824055266496532738528}\\
42 & \texttt{18354203995047060969138151} & 134 & \texttt{58557979706873929433788355}\\
43 & \texttt{18791208275793843875500471} & 135 & \texttt{58995607690467635728659992}\\
44 & \texttt{19228211802839359107199234} & 136 & \texttt{59433424829877633064559607}\\
45 & \texttt{19665215329884874338963873} & 137 & \texttt{59869561653701518718367030}\\
46 & \texttt{20102219610631657245262565} & 138 & \texttt{60305797110366232997864327}\\
47 & \texttt{20539224957281001477325225} & 139 & \texttt{60743904527604531437982051}\\
48 & \texttt{20976231057642937207855461} & 140 & \texttt{61181646917798666743369506}\\
49 & \texttt{21413237158004872938452134} & 141 & \texttt{61617310283583211551650841}\\
50 & \texttt{21850242504654217170445447} & 142 & \texttt{62053600036703618436616865}\\
51 & \texttt{22287246785401000076817631} & 143 & \texttt{62492171559137083121315049}\\
52 & \texttt{22724250312446515308505772} & 144 & \texttt{62929913072061079282226799}\\
53 & \texttt{23161253839492030540281205} & 145 & \texttt{63365103110969803685797358}\\
54 & \texttt{23598258120238813446569563} & 146 & \texttt{63801321329415158997876992}\\
55 & \texttt{24035263466888157678641669} & 147 & \texttt{64240390599124772961811593}\\
56 & \texttt{24472269567250093409163437} & 148 & \texttt{64678313724443055885984508}\\
57 & \texttt{24909275667612029139767798} & 149 & \texttt{65112963222619698661797208}\\
58 & \texttt{25346281014261373371753906} & 150 & \texttt{65548873135968160393810057}\\
59 & \texttt{25783285295008156278133015} & 151 & \texttt{65988475539929668814300182}\\
60 & \texttt{26220288822053671509814470} & 152 & \texttt{66426905757524739257797707}\\
61 & \texttt{26657292349099186741596323} & 153 & \texttt{66861060704171594971763983}\\
62 & \texttt{27094296629845969647878447} & 154 & \texttt{67296201199842820020428144}\\
63 & \texttt{27531301976495313879956896} & 155 & \texttt{67736187996128306439609416}\\
64 & \texttt{27968308076857249610471730} & 156 & \texttt{68175728410677650798922809}\\
65 & \texttt{28405314177219185341084124} & 157 & \texttt{68609759460749839689352310}\\
66 & \texttt{28842319523868529573060775} & 158 & \texttt{69043379008257842365615988}\\
67 & \texttt{29279323804615312479450781} & 159 & \texttt{69482986746056730590840932}\\
68 & \texttt{29716327331660827711120121} & 160 & \texttt{69924368460981353692880012}\\
69 & \texttt{30153330858706342942915226} & 161 & \texttt{70359809876064936293691560}\\
70 & \texttt{30590335139453125849181975} & 162 & \texttt{70791357703932221911489324}\\
71 & \texttt{31027340486102470081282112} & 163 & \texttt{71228204320764617619254789}\\
72 & \texttt{31464346586464405811756285} & 164 & \texttt{71670913623885815406719056}\\
73 & \texttt{31901352686826341542463224} & 165 & \texttt{72111040933952004504815656}\\
74 & \texttt{32338358033475685774197284} & 166 & \texttt{72543493800355947256820413}\\
75 & \texttt{32775362314222468681231190} & 167 & \texttt{72974304146043560753883028}\\
76 & \texttt{33212365841267983911173918} & 168 & \texttt{73411890941237794316905098}\\
77 & \texttt{33649369368313499147610794} & 169 & \texttt{73855465532308049007332176}\\
78 & \texttt{34086373649060282041399607} & 170 & \texttt{74298274427783477103046473}\\
79 & \texttt{34523378995709626307004295} & 171 & \texttt{74735029324019679916754301}\\
80 & \texttt{34960385096071561947666702} & 172 & \texttt{75165705101235334619162984}\\
81 & \texttt{35397391196433497918637463} & 173 & \texttt{75594349576781660415379196}\\
82 & \texttt{35834396543082841509072289} & 174 & \texttt{76023595801279923390844404}\\
83 & \texttt{36271400823829626123566761} & 175 & \texttt{76450522178175231575172444}\\
84 & \texttt{36708404350875136819779650} & 176 & \texttt{76864135612595308625040554}\\
85 & \texttt{37145407877920664058491671} & 177 & \texttt{77246302336992486805167394}\\
86 & \texttt{37582412158667415282166754} & 178 & \texttt{77613257663860083354828908}\\
87 & \texttt{38019417505316842815070159} & 179 & \texttt{78014417926101870463464214}\\
88 & \texttt{38456423605678560391561803} & 180 & \texttt{78443329500014990099443352}\\
89 & \texttt{38893429706041064968206313} & 181 & \texttt{78879413590662700259248153}\\
90 & \texttt{39330435052688932944463707} & 182 & \texttt{90752778779296176935783758}\\
91 & \texttt{39767439333439527085088875} & 183 & \texttt{97600904878812683018771462}\\
92 & \texttt{40204442860475258921208254} &  & \\
\end{longtable}
\end{center}

\section{Correspondence with the formalization}\label{app:lean}

Table~\ref{tab:lean} lists numbered statements of Sections~\ref{sec:intro}--\ref{sec:hostB} together with the corresponding declarations of the accompanying Lean~4 formalization and the files containing them. Names of declarations in the namespace \texttt{Hilbert16} omit the prefix \texttt{Hilbert16.} and the other names are written completely, file paths are relative to the root directory of the Lean project and omit the suffix \texttt{.lean}, and indices in the formalization start at $0$, so that row $i$ of Table~\ref{tab:pairs} and row $k$ of Table~\ref{tab:brackets} are the entries $i-1$ of \nolinkurl{pairRaw} and $k-1$ of \nolinkurl{bracketRaw}.

\begin{center}
\footnotesize
\begin{longtable}{@{}p{0.17\textwidth}p{0.39\textwidth}p{0.38\textwidth}@{}}
\caption{Numbered statements and Lean declarations.}\label{tab:lean}\\
\hline
Statement & Lean declaration & File\\
\hline
\endfirsthead
\hline
Statement & Lean declaration & File\\
\hline
\endhead
\hline
\endfoot
Definition~\ref{defn:cycle} & \nolinkurl{Hilbert16Target.PeriodicSolution},\newline \nolinkurl{Hilbert16Target.PeriodicOrbit},\newline \nolinkurl{Hilbert16Target.IsolatedOrbit},\newline \nolinkurl{Hilbert16Target.HyperbolicLimitCycle} & \nolinkurl{Hilbert16/Target}\\
Theorem~\ref{thm:main} & \nolinkurl{Paper.finite_witnesses} & \nolinkurl{Hilbert16/Paper/FiniteWitnesses}\\
Corollary~\ref{cor:asymptotic} & \nolinkurl{hilbert16_power_lower_bound},\newline \nolinkurl{Hilbert16Target.MainClaim} & \nolinkurl{Hilbert16/Main}\newline \nolinkurl{Hilbert16/Target}\\
Definition~\ref{defn:word} & \nolinkurl{PowerSegment},\newline \nolinkurl{PositiveWord} & \nolinkurl{Paper/PowerSegment}\newline \nolinkurl{Paper/PositiveWord}\\
Lemma~\ref{lem:legal} & \nolinkurl{PositiveWord.word_legal_domain},\newline \nolinkurl{PositiveWord.exists_open_legal_neighborhood},\newline \nolinkurl{PositiveWord.inverse_exponent_product} & \nolinkurl{Paper/PositiveWord/LegalDomain}\newline \nolinkurl{Paper/PositiveWord/OpenNeighborhood}\newline \nolinkurl{Paper/PositiveWord/InverseExponentProduct}\\
Lemma~\ref{lem:stability} & \nolinkurl{finiteSimpleZeroStability} & \nolinkurl{Paper/Analysis/FiniteSimpleZeroStability}\\
Definition~\ref{defn:family} & \nolinkurl{AnalyticFamily} & \nolinkurl{Paper/AnalyticFamily}\\
Theorem~\ref{thm:cert} & \nolinkurl{root_certification} & \nolinkurl{Hilbert16/Paper/RootCertification}\\
Theorem~\ref{thm:compiler} & \nolinkurl{compiler} & \nolinkurl{Paper/Compiler}\\
Proposition~\ref{prop:terminal} & \nolinkurl{terminal_closure} & \nolinkurl{Paper/TerminalClosure}\\
Theorem~\ref{thm:host} & \nolinkurl{squareHost} & \nolinkurl{Paper/Geometry/SquareHost}\\
Lemma~\ref{lem:seed} & \nolinkurl{seed} & \nolinkurl{Paper/Seed}\\
Proposition~\ref{prop:iteration} & \nolinkurl{iteration} & \nolinkurl{Paper/Iteration}\\
Lemma~\ref{lem:perron} & \nolinkurl{perron_estimate} & \nolinkurl{Paper/PerronEstimate}\\
Lemma~\ref{lem:arith} & \nolinkurl{Paper.all_degree_arithmetic} & \nolinkurl{Hilbert16/Paper/AllDegreeArithmetic}\\
Lemma~\ref{lem:exponent} & \nolinkurl{exponent_lower_bound},\newline \nolinkurl{refined_log_margin},\newline \nolinkurl{refined_rational_margin} & \nolinkurl{Hilbert16/Paper/Exponent}\newline \nolinkurl{Hilbert16/Paper/RefinedLogMargin}\newline \nolinkurl{Hilbert16/Paper/RefinedRationalMargin}\\
Tables~\ref{tab:pairs}, \ref{tab:brackets} & \nolinkurl{Certificate.WitnessData.aRaw},\newline \nolinkurl{Certificate.WitnessData.gammaRaw},\newline \nolinkurl{Certificate.WitnessData.pairRaw},\newline \nolinkurl{Certificate.WitnessData.bracketRaw} & \nolinkurl{Paper/Certificate/WitnessData}\\
Lemma~\ref{lem:factor} & \nolinkurl{Certificate.factored_template_derivative},\newline \nolinkurl{template_analytic_odd} & \nolinkurl{Paper/Certificate/FactoredTemplateDerivative}\newline \nolinkurl{Hilbert16/Paper/TemplateAnalyticOdd}\\
Lemma~\ref{lem:oddbracket} & \nolinkurl{odd_bracket_roots} & \nolinkurl{Hilbert16/Analysis/OddBracketRoots}\\
Construction~\ref{cons:interval} & \nolinkurl{Certificate.IntervalProgram.mul},\newline \nolinkurl{Certificate.IntervalProgram.reciprocal},\newline \nolinkurl{Certificate.IntervalProgram.expAt},\newline \nolinkurl{Certificate.IntervalProgram.expInterval} & \nolinkurl{Paper/Certificate/IntervalProgram}\\
Lemma~\ref{lem:sound} & \nolinkurl{Certificate.IntervalProgram.correct} & \nolinkurl{Paper/Certificate/IntervalExponentials}\\
Lemma~\ref{lem:enclosure} & \nolinkurl{Certificate.TemplateProgram.template_enclosure} & \nolinkurl{Paper/Certificate/TemplateEnclosure}\\
Proposition~\ref{prop:checks} & \nolinkurl{Certificate.certificate_checks},\newline \nolinkurl{Certificate.finite_acceptance} & \nolinkurl{Paper/Certificate/CertificateChecks}\newline \nolinkurl{Paper/Certificate/FiniteAcceptance}\\
Lemma~\ref{lem:phase} & \nolinkurl{PhaseReplacement.phase_replacement} & \nolinkurl{Paper/PhaseReplacement}\\
Lemma~\ref{lem:rootmult} & \nolinkurl{rootMultiplication} & \nolinkurl{Paper/RootMultiplication}\\
Lemma~\ref{lem:weight} & \nolinkurl{WeightCoordinate.weight_coordinate} & \nolinkurl{Paper/WeightCoordinate}\\
Lemma~\ref{lem:powerflow} & \nolinkurl{ExactPowerFlows.exact_power_flows} & \nolinkurl{Paper/ExactPowerFlows}\\
Lemma~\ref{lem:transport} & \nolinkurl{frozen_source_transport} & \nolinkurl{Paper/FrozenSourceTransport}\\
Lemma~\ref{lem:pathphase} & \nolinkurl{RegisteredPath.registered_path},\newline \nolinkurl{RegisteredPath.registered_exponential_sum} & \nolinkurl{Paper/RegisteredPath}\newline \nolinkurl{Paper/RegisteredPathExponentialSum}\\
Lemma~\ref{lem:pathfield} & \nolinkurl{RegisteredPath.path_dynamics} & \nolinkurl{Paper/RegisteredPathDynamics}\\
Lemma~\ref{lem:outerweight}, Proposition~\ref{prop:approx} & \nolinkurl{OuterTransport.outer_transport},\newline \nolinkurl{frozen_template_realization} & \nolinkurl{Paper/OuterTransport}\newline \nolinkurl{Paper/FrozenTemplateRealization}\\
Construction~\ref{cons:interleaved}, Lemma~\ref{lem:chainvar} & \nolinkurl{ChainVariation.interleavedPrefix},\newline \nolinkurl{ChainVariation.chain_variation} & \nolinkurl{Paper/ChainVariation}\\
Lemma~\ref{lem:budget} & \nolinkurl{FoldingBudget.folding_budget} & \nolinkurl{Paper/FoldingBudget}\\
Lemma~\ref{lem:height} & \nolinkurl{ChebyshevData.heightControl} & \nolinkurl{Paper/Geometry/HeightControl}\\
Proposition~\ref{prop:polycycle} & \nolinkurl{ChebyshevGeometry.exists_geometry},\newline \nolinkurl{HostSourceFamily.exists_labeling} & \nolinkurl{Paper/Geometry/ChebyshevGeometry}\newline \nolinkurl{Paper/Geometry/SourceLabeling}\\
Lemma~\ref{lem:interp}, Constructions~\ref{cons:field}, \ref{cons:pert} & \nolinkurl{positiveInterpolation},\newline \nolinkurl{HostSourceFamily.polynomialField} & \nolinkurl{Paper/Geometry/PositiveInterpolation}\newline \nolinkurl{Paper/Geometry/HostSourceFamily}\\
Proposition~\ref{prop:saddles} & \nolinkurl{Perturbations.perturbations},\newline \nolinkurl{HostSourceFamily.hostSaddleConnections} & \nolinkurl{Paper/Geometry/Perturbations}\newline \nolinkurl{Paper/Geometry/HostSaddleConnections}\\
Lemma~\ref{lem:samples} & \nolinkurl{EdgeSampling.exists_edgeSamples} & \nolinkurl{Paper/Geometry/EdgeSamples}\\
Lemma~\ref{lem:gapreg} & \nolinkurl{sectionGapRegularity},\newline \nolinkurl{HostSourceFamily.hostLocalSections} & \nolinkurl{Paper/Geometry/SectionGapRegularity}\newline \nolinkurl{Paper/Geometry/HostLocalSections}\\
Lemma~\ref{lem:concentration} & \nolinkurl{quadraticZeroConcentration} & \nolinkurl{Paper/Analysis/Concentration}\\
Proposition~\ref{prop:fullrank} & \nolinkurl{fullRank} & \nolinkurl{Paper/Geometry/FullRank}\\
Lemma~\ref{lem:dulac} & \nolinkurl{factorized_dulac_existence},\newline \nolinkurl{factorized_dulac_scalar_transport} & \nolinkurl{Hilbert16/Paper/FactorizedDulacExistence}\newline \nolinkurl{Hilbert16/Paper/FactorizedDulacTransport}\\
Proposition~\ref{prop:passages} & \nolinkurl{compatibleHostPassages} & \nolinkurl{Paper/Geometry/CompatibleHostPassages}\\
Lemma~\ref{lem:calibration} & \nolinkurl{approximate_calibration} & \nolinkurl{Paper/ApproximateCalibration}\\
Lemma~\ref{lem:gapinverse} & \nolinkurl{parameter_gap_inverse} & \nolinkurl{Hilbert16/Analysis/ParameterGapInverse}\\
Proposition~\ref{prop:return} & \nolinkurl{simultaneousReturn} & \nolinkurl{Paper/Geometry/SimultaneousReturn}\\
Lemma~\ref{lem:persist} & \nolinkurl{hyperbolicPersistence} & \nolinkurl{Paper/Geometry/HyperbolicPersistence}\\
Lemma~\ref{lem:multiplier} & \nolinkurl{returnMultiplierIdentity} & \nolinkurl{Paper/Geometry/ReturnMultiplierIdentity}\\
Lemma~\ref{lem:isolation} & \nolinkurl{floquetOrbitBridge},\newline \nolinkurl{returnOrbitIsolation},\newline \nolinkurl{c1FirstReturnCycle},\newline \nolinkurl{closedReturnImages} & \nolinkurl{Paper/Geometry/FloquetOrbitBridge}\newline \nolinkurl{Paper/Geometry/ReturnOrbitIsolation}\newline \nolinkurl{Paper/Geometry/C1FirstReturnCycle}\newline \nolinkurl{Paper/Geometry/ClosedReturnImages}\\
\end{longtable}
\end{center}

\end{document}